\documentclass[11pt,draft]{amsart}

\usepackage{a4,amssymb,cite}

\usepackage{pstricks,pst-node,pst-text}

\usepackage{graphicx,gastex}

\usepackage[]{color}

\DeclareMathOperator{\Con}{Con}

\DeclareMathOperator{\Eq}{Eq}

\DeclareMathOperator{\Perm}{Perm}

\DeclareMathOperator{\Stab}{Stab}

\DeclareMathOperator{\Sub}{Sub}

\DeclareMathOperator{\var}{var}

\newtheorem{theorem}{Theorem}[section]

\newtheorem{proposition}[theorem]{Proposition}

\newtheorem{lemma}[theorem]{Lemma}

\newtheorem{corollary}[theorem]{Corollary}

\newtheorem{remark}[theorem]{Remark}

\newtheorem{observation}[theorem]{Observation}

\theoremstyle{definition}

\newtheorem{question}[theorem]{Question}

\newcommand{\contsection}[2]{\small\ref{#2}.\enskip#1\dotfill\pageref{#2}}

\newcommand{\contsubsection}[2]{\footnotesize\ref{#2}.~#1\ (\pageref{#2}).}

\newcommand{\contstarsection}[2]{\small#1\dotfill\pageref{#2}}

\makeatletter

\renewcommand*\subjclass[2][2010]{\def\@subjclass{#2}\@ifundefined{subjclassname@#1}{\ClassWarning{\@classname}{Unknown edition (#1) of Mathematics Subject Classification; using '2010'.}}{\@xp\let\@xp\subjclassname\csname subjclassname@#1\endcsname}}

\renewcommand{\subjclassname}{\textup{2010} Mathematics Subject Classification}

\makeatother

\begin{document}

\title[Special elements in lattices of semigroup varieties]{Special elements in lattices\\
of semigroup varieties}

\thanks{At various stages of writing this article, the work was supported by the Ministry of Education and Science of the Russian Federation (projects 2248/2014 and 1.6018.2017/8.9), by grant of the President of the Russian Federation for supporting of leading scientific schools of the Russian Federation (project 5161.2014.1) and by Russian Foundation for Basic Research (grants 14-01-00524 and 17-01-00551).}

\author{B.\,M.\,Vernikov}

\address{Ural Federal University, Institute of Natural Sciences and Mathematics, Lenina 51, 620000 Ekaterinburg, Russia}

\email{bvernikov@gmail.com}

\date{}

\begin{abstract}
We survey results concerning special elements of nine types (modular, lower-modular, upper-modular, cancellable, distributive, codistributive, standard, costandard and neutral elements) in the lattice of all semigroup varieties and certain its sublattices, mainly in the lattices of all commutative varieties and of all overcommutative ones.
\end{abstract}



\maketitle

{\centering\section*{Contents}}

\begin{itemize}
\item[]\contsection{Introduction}{section intr}
\item[]\contsection{Preliminary results}{section prel}
\begin{quote}
\contsubsection{$I$-elements and $Q$-elements of lattices}{subsection abstr,SEM,Com:I-,Q-elem} \contsubsection{Modular and upper-modular elements in lattices of equivalence relations}{subsection Eq:mod,umod} \contsubsection{Special elements in subgroup lattices of finite symmetric groups}{subsection Sub:all-but-umod,lmod} \contsubsection{Special elements in congruence lattices of $G$-sets}{subsection gs:distr,codistr,canc,stand,costand,neutr} \contsubsection{Upper-modular and codistributive elements: interrelations between lattice identities and a hereditary property}{subsection umod,codistr:ident&hered}
\end{quote}
\item[]\contsection{The lattice $\mathbb{SEM}$}{section SEM}
\begin{quote}
\contsubsection{$I$-varieties}{subsection SEM:I-varieties} \contsubsection{Lower-modular varieties}{subsection SEM:lmod} \contsubsection{Distributive and standard varieties}{subsection SEM:distr,stand} \contsubsection{Costandard and neutral varieties}{subsection SEM:costand,neutr} \contsubsection{An application to definable varieties}{subsection SEM:definability} \contsubsection{Modular varieties}{subsection SEM:cmod} \contsubsection{Cancellable varieties}{subsection SEM:canc} \contsubsection{Upper-modular varieties}{subsection SEM:umod} \contsubsection{Varieties that are both modular and upper-modular}{subsection SEM:cumod} \contsubsection{Codistributive varieties}{subsection SEM:codistr}
\end{quote}
\item[]\contsection{The lattice $\mathbb{COM}$}{section COM}
\begin{quote}
\contsubsection{$\mathbb{COM}$-lower-modular varieties}{subsection COM:lmod} \contsubsection{$\mathbb{COM}$-distributive and $\mathbb{COM}$-standard varieties}{subsection COM:distr,stand} \contsubsection{$\mathbb{COM}$-neutral varieties}{subsection COM:neutr} \contsubsection{$\mathbb{COM}$-modular varieties}{subsection COM:cmod} \contsubsection{$\mathbb{COM}$-cancellable varieties}{subsection COM:canc} \contsubsection{$\mathbb{COM}$-upper-modular and $\mathbb{COM}$-codistributive varieties}{subsection COM:umod,codistr} \contsubsection{$\mathbb{COM}$-co\-standard varieties}{subsection COM:costand}
\end{quote}
\item[]\contsection{Lattices located between $\mathbb{SEM}$ and $\mathbb{COM}$}{section between SEM and COM}
\begin{quote}
\contsubsection{Subvariety lattices of overcommutative varieties}{subsection between-SEM-and-COM:L(oc)} \contsubsection{The lattice $\mathbb{PERM}$}{subsection between-SEM-and-COM:PERM}
\end{quote}
\item[]\contsection{The lattice $\mathbb{OC}$}{section OC}
\item[]\contstarsection{Acknowledgments}{thanks}
\item[]\contstarsection{\refname}{bibl}
\end{itemize}

\section{Introduction}
\label{section intr}

This work is an extended version of the survey~\cite{Vernikov-15}. The work is regularly updated and modified as new results and/or articles appear.

The collection of all semigroup varieties forms a lattice with respect to class-theoretical inclusion. This lattice will be denoted by $\mathbb{SEM}$. The lattice $\mathbb{SEM}$ has been intensively studied since the beginning of 1960s. A systematic overview of the material accumulated here is given in the survey~\cite{Shevrin-Vernikov-Volkov-09}.

The lattice $\mathbb{SEM}$ has an extremely complicated structure. In particular, it contains an anti-isomorphic copy of the partition lattice over a countably infinite set~\cite{Burris-Nelson-71-infinite,Jezek-76}, and therefore, does not satisfy any non-trivial lattice identity. Identities in subvariety lattices of semigroup varieties were intensively examined in many articles. These articles contain a number of interesting and deep results (see~\cite[Section~11]{Shevrin-Vernikov-Volkov-09}). The next natural step is to consider varieties that guarantee, so to speak, `nice lattice behavior' in their neighborhood. Specifically, our attention is to study special elements of different types in the lattice $\mathbb{SEM}$.

We will consider nine types of special elements: modular, lower-modular, upper-modular, cancellable, distributive, codistributive, standard, costandard and neutral elements. Recall the corresponding definitions. An element $x$ of a lattice $\langle L;\vee,\wedge\rangle$ is called
\begin{align*}
&\text{\emph{modular} if}\quad&&\forall\,y,z\in L\colon\quad y\le z\longrightarrow(x\vee y)\wedge z=(x\wedge z)\vee y;\\
&\text{\emph{lower-modular} if}\quad&&\forall\,y,z\in L\colon\quad x\le y\longrightarrow x\vee(y\wedge z)=y\wedge(x\vee z);\\
&\text{\emph{cancellable} if}\quad&&\forall\,y,z\in L\colon\quad x\vee y=x\vee z\ \&\ x\wedge y=x\wedge z\longrightarrow y=z;\\
&\text{\emph{distributive} if}\quad&&\forall\,y,z\in L\colon\quad x\vee(y\wedge z)=(x\vee y)\wedge(x\vee z);\\
&\text{\emph{standard} if}\quad&&\forall\,y,z\in L\colon\quad(x\vee y)\wedge z=(x\wedge z)\vee(y\wedge z);
\end{align*}
\emph{neutral} if, for all $y,z\in L$, the sublattice of $L$ generated by $x$, $y$ and $z$ is distributive. It is well known (see~\cite[Theorem~254 on p.\,226]{Gratzer-11}, for instance) that an element $x\in L$ is neutral if and only if
$$
\forall\,y,z\in L\colon\quad(x\vee y)\wedge(y\vee z)\wedge(z\vee x)=(x\wedge y)\vee(y\wedge z)\vee(z\wedge x).
$$
\emph{Upper-modular, codistributive} and \emph{costandard} elements are defined dually to lower-modular, distributive and standard ones respectively.

Special elements play an important role in the general lattice theory (see~\cite[Section~III.2]{Gratzer-11}, for instance). In particular, it is well known that if $a$ is a neutral element in a lattice $L$ then $L$ is decomposable into a subdirect product of the principal ideal and the principal filter of $L$ generated by $a$ (see~\cite[Theorem~254 on p.\,226]{Gratzer-11}, for instance). Thus, the knowledge of which elements of a lattice are neutral gives essential information on the structure of the lattice as a whole. A valuable information about special elements of types mentioned above may be found in~\cite{Gratzer-11,Seselja-Tepavcevic-01,Stern-99}.

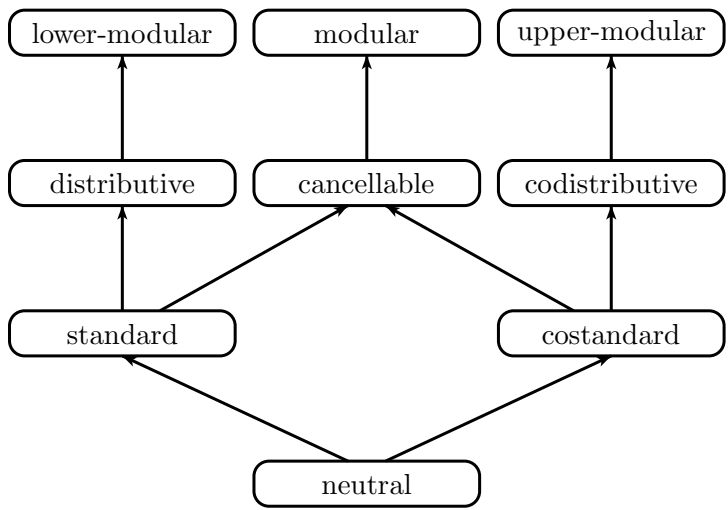
\begin{figure}[tbh]
\begin{center}
\unitlength=1mm
\linethickness{0.4pt}
\begin{picture}(95,66)
\gasset{linewidth=0.4}
\drawoval(15,23,30,6,2)
\drawoval(15,43,30,6,2)
\drawoval(15,63,30,6,2)
\drawoval(47.5,3,30,6,2)
\drawoval(47.5,43,30,6,2)
\drawoval(47.5,63,30,6,2)
\drawoval(80,23,30,6,2)
\drawoval(80,43,30,6,2)
\drawoval(80,63,30,6,2)
\put(15,23){\makebox(0,0)[cc]{standard}}
\put(15,43){\makebox(0,0)[cc]{distributive}}
\put(15,63){\makebox(0,0)[cc]{lower-modular}}
\put(47.5,3){\makebox(0,0)[cc]{neutral}}
\put(47.5,43){\makebox(0,0)[cc]{cancellable}}
\put(47.5,63){\makebox(0,0)[cc]{modular}}
\put(80,23){\makebox(0,0)[cc]{costandard}}
\put(80,43){\makebox(0,0)[cc]{codistributive}}
\put(80,63){\makebox(0,0)[cc]{upper-modular}}
\gasset{AHnb=1,AHLength=2}
\drawline(45,6)(15,20)
\drawline(50,6)(80,20)
\drawline(15,26)(15,40)
\drawline(20,26)(45,40)
\drawline(75,26)(50,40)
\drawline(80,26)(80,40)
\drawline(15,46)(15,60)
\drawline(47.5,46)(47.5,60)
\drawline(80,46)(80,60)
\end{picture}
\caption{Interrelations between types of elements in abstract lattices}
\label{relations abstract}
\end{center}
\end{figure}

There is a number of interrelations between types of elements we consider. It is evident that a neutral element is both standard and costandard; a standard or costandard element is cancellable; a cancellable element is modular; a [co]distributive element is lower-modular [upper-modular]. It is well known also that a [co]standard element is [co]distributive (see~\cite[Theorem~253 on p.\,224]{Gratzer-11}, for instance). These interrelations between types of elements in abstract lattices are shown in Fig.~\ref{relations abstract}.

In fact, first information about special elements in the lattice $\mathbb{SEM}$ were obtained (in the implicit form) in the articles~\cite{Aizenstat-74,Melnik-69,Melnik-70,Melnik-73,Salij-69}. There were no any explicit references to special elements in these articles. But it immediately follows from the results obtained there that particular concrete varieties are special elements of some or another type in $\mathbb{SEM}$. Translated into the language of special elements, the results of these works can be formulated as follows. It is proved independently in~\cite{Melnik-69} and~\cite{Salij-69} that the variety of semilattices $\mathbf{SL}$ is a modular element of $\mathbb{SEM}$. It is verified in~\cite{Melnik-70} that $\mathbf{SL}$ is distributive in $\mathbb{SEM}$. The fact that the variety $\mathbf{ZM}$ of semigroups with zero multiplication is modular and distributive in $\mathbb{SEM}$ is checked in~\cite{Melnik-73}. Finally, it is verified in~\cite{Aizenstat-74} that the lattice $\mathbb{SEM}$ is 0-\emph{disributive}, i.e., satisfies the implication
$$
x\wedge y=x\wedge z=0\longrightarrow x\wedge(y\vee z)=0.
$$
This immediately implies that all atoms of the lattice $\mathbb{SEM}$ (in particular, the varieties $\mathbf{SL}$ and $\mathbf{ZM}$) are codistributive in $\mathbb{SEM}$.

First explicit results about special elements in $\mathbb{SEM}$ were obtained in the articles~\cite{Jezek-McKenzie-93,Vernikov-Volkov-88}. In~\cite{Vernikov-Volkov-88}, certain sufficient condition for a semigroup variety to be modular or lower-modular element of $\mathbb{SEM}$ is found. In~\cite{Jezek-McKenzie-93}, Je\v{z}ek and McKenzie examined modular elements in $\mathbb{SEM}$\footnote{Note that the paper~\cite{Jezek-McKenzie-93} deals with the lattice of equational theories of semigroups, i.e., the dual of $\mathbb{SEM}$ rather than the lattice $\mathbb{SEM}$ itself. When reproducing results from~\cite{Jezek-McKenzie-93}, we adapt them to the terminology of the present article.}. They rediscover the sufficient condition for modular elements of $\mathbb{SEM}$ mentioned in~\cite{Vernikov-Volkov-88} and find a strong necessary condition for a semigroup variety to be modular element of $\mathbb{SEM}$. Note that the mentioned results of~\cite{Jezek-McKenzie-93,Vernikov-Volkov-88} play an auxiliary role in these works. 
 
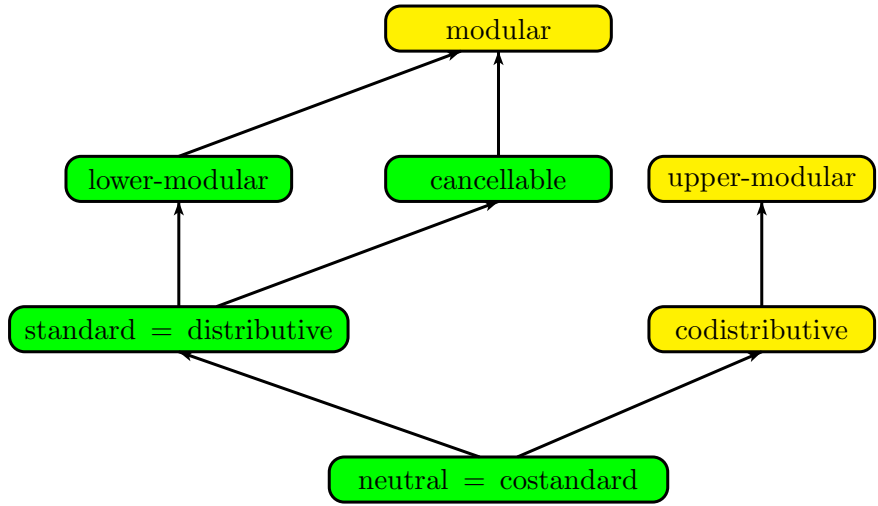
\begin{figure}[tbh]
\begin{center}
\unitlength=1mm
\linethickness{0.4pt}
\begin{picture}(115,66)
\gasset{linewidth=0.4}
\drawoval[fillcolor=green](22.5,23,45,6,2)
\drawoval[fillcolor=green](22.5,43,30,6,2)
\drawoval[fillcolor=green](65,3,45,6,2)
\drawoval[fillcolor=green](65,43,30,6,2)
\drawoval[fillcolor=yellow](65,63,30,6,2)
\drawoval[fillcolor=yellow](100,23,30,6,2)
\drawoval[fillcolor=yellow](100,43,30,6,2)
\put(22.5,23){\makebox(0,0)[cc]{standard\enskip=\enskip distributive}}
\put(22.5,43){\makebox(0,0)[cc]{lower-modular}}
\put(65,3){\makebox(0,0)[cc]{neutral\enskip=\enskip costandard}}
\put(65,43){\makebox(0,0)[cc]{cancellable}}
\put(65,63){\makebox(0,0)[cc]{modular}}
\put(100,23){\makebox(0,0)[cc]{codistributive}}
\put(100,43){\makebox(0,0)[cc]{upper-modular}}
\gasset{AHnb=1,AHLength=2}
\drawline(62.5,6)(22.5,20)
\drawline(67.5,6)(100,20)
\drawline(22.5,26)(22.5,40)
\drawline(27.5,26)(65,40)
\drawline(100,26)(100,40)
\drawline(22.5,46)(60,60)
\drawline(65,46)(65,60)
\end{picture}
\caption{Interrelations between types of elements in $\mathbb{SEM}$}
\label{relations SEM}
\end{center}
\end{figure}

A systematic examination of special elements in $\mathbb{SEM}$ is the objective of the articles~\cite{Gusev-Skokov-Vernikov-18,Shaprynskii-12-mod-lmod,Shaprynskii-12-I-elem,Shaprynskii-Skokov-Vernikov-19,Shaprynskii-Vernikov-10,Skokov-Vernikov-19,Vernikov-07-cmod,Vernikov-07-lmod,Vernikov-08-lmod,Vernikov-08-umod1,Vernikov-08-umod2,Vernikov-11,Vernikov-Shaprynskii-10,Vernikov-Volkov-06,Volkov-05}; see also~\cite[Section~14]{Shevrin-Vernikov-Volkov-09}. In short, the mentioned articles contain complete descriptions of lower-modular, cancellable, distributive, standard, costandard and neutral elements of the lattice $\mathbb{SEM}$\footnote{To prevent a possible confusion, we note that the description of standard elements of $\mathbb{SEM}$ is not formulated explicitly anywhere but readily follows from results of~\cite{Vernikov-Shaprynskii-10}, see a comment to Theorem~\ref{SEM distr,stand} below.} and essential information (such as strong necessary conditions and descriptions in wide and important partial cases) about modular, upper-modular and codistributive elements of this lattice. In particular, it turns out that there are some interrelations between special elements of different types in $\mathbb{SEM}$ that do not hold in abstract lattices. Namely, an element of $\mathbb{SEM}$ is standard if and only if it is distributive; is costandard if and only if it is neutral; is modular whenever it is lower-modular. Interrelations between types of elements in the lattice $\mathbb{SEM}$ are shown in Fig.~\ref{relations SEM}. All other possible interrelations between types of elements under consideration are not the case. Corresponding examples may be easily extracted from results formulated below. In Figures~\ref{relations SEM} and~\ref{relations COM} we color the oval corresponding to the type of elements in green color if elements of this type in the corresponding lattice are completely determined, and in yellow color, if an essential information about such elements is known.

The lattice $\mathbb{SEM}$ contains a number of wide and important sublattices (see~\cite[Section~1 and Chapter~2]{Shevrin-Vernikov-Volkov-09}). It is natural to examine special elements in these sublattices. One of the most important sublattices of $\mathbb{SEM}$ is the lattice $\mathbb{COM}$ of all commutative semigroup varieties. It follows from results of~\cite{Burris-Nelson-71-finite} that this lattice contains an isomorphic copy of any finite lattice, and therefore, does not satisfy any non-trivial lattice identity. On the other hand, the lattice $\mathbb{COM}$ is known to be countably infinite~\cite{Perkins-69} and can be characterized~\cite{Kisielewicz-94} (see also~\cite[Section~8]{Shevrin-Vernikov-Volkov-09}). Special elements in the lattice $\mathbb{COM}$ are examined in~\cite{Shaprynskii-11,Shaprynskii-12-mod-lmod,Vernikov-17} where lower-modular, upper-modular, distributive, codistributive, standard, costandard and neutral elements of $\mathbb{COM}$ are completely determined and an essential information about modular elements of this lattice is obtained. As in the case of the lattice $\mathbb{SEM}$, it turns out that an element of $\mathbb{COM}$ is standard if and only if it is distributive, and is modular whenever it is lower-modular. This deep analogy between properties of special elements in the lattices $\mathbb{SEM}$ and $\mathbb{COM}$ is not accidental. As we will seen below, these similar results easily follow from some general fact concerning properties of special elements in lattices of subvarieties of overcommutative varieties (evidently, both the lattices $\mathbb{SEM}$ and $\mathbb{COM}$ are subvariety lattices of two `extremal' overcommutative varieties, namely the variety of all semigroups and the variety of all commutative semigroups respectively). But this analogy does not extend to all types of special elements. In contrast with the case of $\mathbb{SEM}$, it turns out that the properties of being neutral and costandard elements of the lattice $\mathbb{COM}$ are not equivalent, whereas the properties of being codistributive and upper-modular elements of this lattice are, on the contrary, equivalent. Interrelations between types of elements in the lattice $\mathbb{COM}$ are shown in Fig.~\ref{relations COM}. We do not mention cancellable elements in this figure because there no any information about these elements in the lattice $\mathbb{COM}$ so far (except necessary conditions for modular elements of the lattice $\mathbb{COM}$ that evidently are also necessary conditions for cancellable elements of this lattice). No interrelations between types of elements in $\mathbb{COM}$ not specified in Fig.~\ref{relations COM} hold. As in the case of the lattice $\mathbb{SEM}$, corresponding examples may be easily extracted from results formulated below.{\sloppy

}

\begin{figure}[tbh]
\begin{center}
\unitlength=1mm
\linethickness{0.4pt}
\begin{picture}(115,86)
\gasset{AHnb=0,linewidth=0.4}
\drawoval[fillcolor=green](22.5,43,45,6,2)
\drawoval[fillcolor=green](22.5,63,30,6,2)
\drawoval[fillcolor=yellow](22.5,83,30,6,2)
\drawoval[fillcolor=green](50,3,30,6,2)
\drawoval[fillcolor=green](50,23,30,6,2)
\drawoval[fillcolor=green](85,43,60,6,2)
\put(22.5,43){\makebox(0,0)[cc]{standard\enskip=\enskip distributive}}
\put(22.5,63){\makebox(0,0)[cc]{lower-modular}}
\put(22.5,83){\makebox(0,0)[cc]{modular}}
\put(50,3){\makebox(0,0)[cc]{neutral}}
\put(50,23){\makebox(0,0)[cc]{costandard}}
\put(85,43){\makebox(0,0)[cc]{codistributive\enskip=\enskip upper-modular}}
\gasset{AHnb=1,AHLength=2}
\drawline(50,6)(50,20)
\drawline(47.5,26)(22.5,40)
\drawline(52.5,26)(85,40)
\drawline(22.5,46)(22.5,60)
\drawline(22.5,66)(22.5,80)
\end{picture}
\caption{Interrelations between types of elements in $\mathbb{COM}$}
\label{relations COM}
\end{center}
\end{figure}

Recall that a semigroup variety is called \emph{permutative} if it satisfies a \emph{permutational} identity, i.e., an identity of the type
$$
x_1x_2\cdots x_n\approx x_{1\pi}x_{2\pi}\cdots x_{n\pi}
$$
where $\pi$ is a non-trivial permutation on the set $\{1,2,\dots,n\}$. This identity will be denoted by $p_n[\pi]$. The number $n$ is called the \emph{length} of this identity. The collection of all permutative varieties forms a sublattice $\mathbb{PERM}$ of the lattice $\mathbb{SEM}$. This lattice is located between $\mathbb{SEM}$ and $\mathbb{COM}$. It seems quite natural to examine special elements in $\mathbb{PERM}$. There are no published results here so far but some information about modular and lower-modular elements in the lattice $\mathbb{PERM}$ is found in PhD thesis by Shaprynski\v{\i}~\cite{Shaprynskii-15}.

The `antipode' of the lattice $\mathbb{COM}$ is the lattice $\mathbb{OC}$ of all \emph{overcommutative} semigroup varieties (i.e., varieties containing the variety of all commutative semigroups). It is well known that the lattice $\mathbb{SEM}$ is the disjoint union of $\mathbb{OC}$ and the lattice of all \emph{periodic} semigroup varieties (i.e., varieties consisting of periodic semigroups). Results of the papers~\cite{Jezek-McKenzie-93,Vernikov-07-lmod,Vernikov-08-umod1} imply that if a semigroup variety $\mathbf V$ is an element of one of the nine types mentioned above in the lattice $\mathbb{SEM}$ and $\mathbf V$ is different from the variety of all semigroups then $\mathbf V$ is a periodic variety (more general fact is proved in~\cite{Shaprynskii-12-I-elem}, see Proposition~\ref{I-varieties are periodic} below). Thus, an examination of special elements of all mentioned types in $\mathbb{SEM}$ a~priori can not give any essential information about the lattice $\mathbb{OC}$. Note that the lattice $\mathbb{OC}$ contains an isomorphic copy of any finite lattice~\cite{Volkov-94}, whence it does not satisfy any non-trivial lattice identity. Overcommutative varieties whose lattice of overcommutative subvarieties satisfies a particular lattice identity were intensively studied (see~\cite[Subsection~5.2]{Shevrin-Vernikov-Volkov-09} and the article~\cite{Shaprynskii-13}). All these arguments make the examination of special elements of $\mathbb{OC}$ very natural. Such an examination has been started in the article~\cite{Vernikov-01}. It is proved there that the properties of being a distributive, a codistributive, a standard, a costandard and a neutral element of the lattice $\mathbb{OC}$ are equivalent, and a certain characterization of corresponding overcommutative varieties is proposed. But this description turns out to be incorrect (while the result that the five mentioned conditions are equivalent is true). The correct description of distributive, codistributive, standard, costandard and neutral elements of the lattice $\mathbb{OC}$ is contained in the article~\cite{Shaprynskii-Vernikov-11}. Cancellable elements of the lattice $\mathbb{OC}$ are classified in~\cite{Shaprynskii-Vernikov-21}. More precisely, it is verified there that a property to be a cancellable element of $\mathbb{OC}$ is equivalent to the property to be a [co]distributive, [co]standard or neural element of this lattice. There are no any information about modular, lower-modular or upper-modular elements of $\mathbb{OC}$ so far.

Note that there is an interesting information about special elements in lattices of varieties of semigroups with different additional unary or nullary operations. So, a number of examples of neutral elements in the lattice of all completely regular semigroup varieties is given in~\cite{Trotter-89}; here completely regular semigroups are considered as \emph{unary semigroups}, that is, semigroups with a naturally defined additional unary operation (see~\cite{Petrich-Reilly-99} or~\cite[Section~6]{Shevrin-Vernikov-Volkov-09}, for instance). Another type of unary semigroups which includes completely regular semigroups as a partial case is epigroups (see~\cite{Shevrin-94,Shevrin-05} or~\cite[Section~2]{Shevrin-Vernikov-Volkov-09}, for instance). Special elements of all mentioned above types in the lattice of epigroup varieties are examined in~\cite{Shaprynskii-Skokov-Vernikov-16,Shaprynskii-Skokov-Vernikov-19,Skokov-15,Skokov-16,Skokov-18}. The works~\cite{Gusev-18,Gusev-20,Gusev-Lee-21+} are devoted to special elements of several types in the lattice of monoid varieties. But all these results are beyond the scope of this survey.

The survey consists of six sections. In Section~\ref{section prel}, we provide some preliminary results about special elements in abstract lattices, lattices of equivalence relations, subgroup lattices of finite symmetric groups, congruence lattices of $G$-sets and the lattices $\mathbb{SEM}$ and $\mathbb{COM}$. These preliminary results play an important role in the proofs of the results that we survey in Sections~\ref{section SEM}--\ref{section OC}. In Sections~\ref{section SEM} and~\ref{section COM}, we overview results about special elements in the lattices $\mathbb{SEM}$ and $\mathbb{COM}$ respectively. Section~\ref{section between SEM and COM} contains results about modular and lower-modular elements in lattices located between $\mathbb{SEM}$ and $\mathbb{COM}$, namely in subvariety lattices of overcommutative varieties and in the lattice $\mathbb{PERM}$. Finally, Section~\ref{section OC} is devoted to special elements in the lattice $\mathbb{OC}$. Sections~\ref{section SEM} and~\ref{section between SEM and COM} contain also several open questions. 

\section{Preliminary results}
\label{section prel}

\subsection{$I$-elements and $Q$-elements of lattices}
\label{subsection abstr,SEM,Com:I-,Q-elem}

Almost all types of special elements introduced above (namely, all of them except cancellable elements) are defined by the following general scheme. We take a particular lattice identity and consider it as an open formula. Then, one of the variables is left free while all the others are subjected to a universal quantifier. As a result, we obtain a first order formula $\Phi(x)$ with one free variable $x$ in the language of lattice operations. An element $w$ of a lattice $L$ is said to be a special element of corresponding type of $L$ if the sentence $\Phi(w)$ is true. It is evident that neutral, [co]standard and [co]distributive elements are defined just by this scheme. The definitions of modular, lower-modular and upper-modular elements may be easily reformulated in the framework of this approach as well because the modular law may be written as an identity.

One can formalize the approach discussed in the previous paragraph. We will write lattice or semigroup terms (rather than letters) in bold and connect two sides of a lattice or semigroup identity by the symbol $\approx$. The symbol $=$ will denote, among other things, the equality relation on a lattice or semigroup. Let $\varepsilon$ be a lattice identity of the form $\mathbf{s\approx t}$ where terms $\mathbf s$ and $\mathbf t$ depend on variables $x_1,\dots,x_n$, and $1\le i\le n$. To write [quasi-]identities more compact, we put $X_n^i=\{x_1,\dots,x_{i-1},x,x_{i+1},\dots,x_n\}$. An element $x$ of a lattice $L$ is called an $(\varepsilon,i)$-\emph{element} of $L$ if
$$
\forall\,x_1,\dots,x_{i-1},x_{i+1},\dots,x_n\in L\colon\quad\mathbf s(X_n^i)=\mathbf t(X_n^i).
$$
An element of a lattice $L$ is called an $I$-\emph{element} of $L$ if it is an $(\varepsilon,i)$-element of $L$ for some non-trivial identity $\varepsilon$ depending on variables $x_1,\dots,x_n$ and some $1\le i\le n$.

For an element $a$ of a lattice $L$, we put $(a]=\{x\in L\mid x\le a\}$. If $a\in L$ and the lattice $(a]$ satisfies the identity $\mathbf p(x_1,\dots,x_n)\approx\mathbf q(x_1,\dots,x_n)$ then
$$
\mathbf p(a\wedge x_1,\dots,a\wedge x_n)=\mathbf q(a\wedge x_1,\dots,a\wedge x_n)
$$
for all $x_1,\dots,x_n\in L$ because $a\wedge x_1,\dots,a\wedge x_n\in(a]$. Therefore, in this situation, $a$ is an $(\varepsilon,n+1)$-element of $L$ with the following $\varepsilon$:
$$
\mathbf p(x_{n+1}\wedge x_1,\dots,x_{n+1}\wedge x_n)\approx\mathbf q(x_{n+1}\wedge x_1,\dots,x_{n+1}\wedge x_n).
$$
So, we have the following

\begin{observation}
\label{identities and I-elements}
If $w$ is an element of a lattice $L$ and the ideal $(w]$ of $L$ satisfies some non-trivial lattice identity then $w$ is an $I$-element of $L$.
\end{observation}

The subvariety lattice of a variety $\mathbf V$ is denoted by $L(\mathbf V)$. A semigroup variety $\mathbf V$ is called an $I$-\emph{variety} if it is an $I$-element of the lattice $\mathbb{SEM}$. The following assertion is a specialization of Observation~\ref{identities and I-elements} for the lattice $\mathbb{SEM}$.

\begin{corollary}
\label{identities and I-varieties}
If $\mathbf V$ is a semigroup variety and the lattice $L(\mathbf V)$ satisfies some non-trivial lattice identity then $\mathbf V$ is an $I$-variety.
\end{corollary}

We will denote by $\var\Sigma$ the semigroup variety given by the identity system $\Sigma$. The converse statement to Corollary~\ref{identities and I-varieties} is not true. Indeed, the variety $\var\{x^2\approx0\}$ is a modular and a lower-modular element of $\mathbb{SEM}$ (see Theorems~\ref{SEM lmod} and~\ref{SEM cmod suf} below) but its subvariety lattice does not satisfy any non-trivial identity~\cite{Jezek-76}.

The following fact turns out to be very helpful.

\begin{lemma}[\!\!{\mdseries\cite[Corollary~2.1]{Shaprynskii-11}}]
\label{join with neutr atom I-elem}
Let $w$ be an atom and a neutral element of a lattice $L$, $\varepsilon$ a lattice identity that holds in the $2$-element lattice and depends on variables $x_1,\dots,x_n$, and $1\le i\le n$. An element $x\in L$ is an $(\varepsilon,i)$-element of $L$ if and only if the element $x\vee w$ has the same property.
\end{lemma}

It is well known that the variety $\mathbf{SL}$ is an atom of the lattice $\mathbb{SEM}$ (see~\cite[Section~1]{Shevrin-Vernikov-Volkov-09}, for instance) and a neutral element of this lattice (see~\cite[Proposition~4.1]{Volkov-05} or Theorem~\ref{SEM costand,neutr} below). In particular, $\mathbf{SL}$ is a neutral atom of $\mathbb{COM}$. Thus, Lemma~\ref{join with neutr atom I-elem} implies the following

\begin{corollary}
\label{join with SL I-elem}
Let $\varepsilon$ be a lattice identity that holds in the $2$-element lattice and depends on variables $x_1,\dots,x_n$, and $1\le i\le n$. A \textup[commutative\textup] semigroup variety $\mathbf V$ is an $(\varepsilon,i)$-element of the lattice $\mathbb{SEM}$ \textup[respectively $\mathbb{COM}$\textup] if and only if the variety $\mathbf{V\vee SL}$ has the same property.
\end{corollary}

Note that a number of partial cases of Lemma~\ref{join with neutr atom I-elem} and Corollary~\ref{join with SL I-elem} for special elements of different concrete types were proved earlier in~\cite{Vernikov-07-lmod,Vernikov-11,Vernikov-Shaprynskii-10,Vernikov-Volkov-06,Volkov-05}.

The notion of $I$-elements does not cover the notion of cancellable elements. It seems natural to generalize the former notion in the following way. Let $\xi$ be a lattice quasi-identity of the form
$$
\mathop\&\limits_{k=1}^m\mathbf s_k\approx\mathbf t_k\longrightarrow \mathbf{s\approx t}
$$
where terms $\mathbf s_1,\dots,\mathbf s_m,\mathbf t_1,\dots,\mathbf t_m,\mathbf s$ and $\mathbf t$ depend on variables $x_1,\dots,x_n$, and $1\le i\le n$. An element $x$ of a lattice $L$ is called a $(\xi,i)$-\emph{element} of $L$ if
$$
\forall\,x_1,\dots,x_{i-1},x_{i+1},\dots,x_n\in L\colon\quad\mathop\&\limits_{k=1}^m\mathbf s_k(X_n^i)=\mathbf t_k(X_n^i)\longrightarrow\mathbf s(X_n^i)=\mathbf t(X_n^i).
$$
An element of a lattice $L$ is called a $Q$-\emph{element} of $L$ if it is a $(\xi,i)$-element of $L$ for some non-trivial quasi-identity $\xi$ depending on variables $x_1,\dots,x_n$ and some $1\le i\le n$. Clearly, cancellable elements are $Q$-elements.

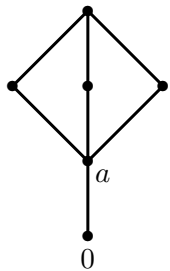
\begin{figure}[tbh]
\begin{center}
\unitlength=1mm
\linethickness{0.4pt}
\begin{picture}(20,36)
\gasset{AHnb=0,linewidth=0.4}
\drawline(10,5)(10,35)(0,25)(10,15)(20,25)(10,35)
\put(0,25){\circle*{1.33}}
\put(10,5){\circle*{1.33}}
\put(10,15){\circle*{1.33}}
\put(10,25){\circle*{1.33}}
\put(10,35){\circle*{1.33}}
\put(20,25){\circle*{1.33}}
\put(10,2){\makebox(0,0)[cc]{0}}
\put(12,13){\makebox(0,0)[cc]{$a$}}
\end{picture}
\caption{The lattice $M_3$ with the new least element adjoined}
\label{M_3 with 0}
\end{center}
\end{figure}

Unfortunately, the analogue of Lemma~\ref{join with neutr atom I-elem} for $Q$-elements is not the case. More precisely, the join of a $(\xi,i)$-element of a lattice $L$ and a neutral atom of $L$ need not to be a $(\xi,i)$-element of $L$. Indeed, let $\xi$ be the following quasi-identity:
$$
x_1\approx x_2\wedge x_3\longrightarrow(x_2\wedge x_3)\vee x_4\approx(x_2\vee x_4)\wedge(x_3\vee x_4).
$$
Further, let $L$ be the 5-element modular non-distributive lattice $M_3$ with the new least element adjoined (see Fig.~\ref{M_3 with 0}). We denote this new least element by~0 and a unique atom of $L$ by $a$. Then it is easy to see that the atom $a$ is neutral,~0 is a $(\xi,1)$-element of $L$ but the element $a=0\vee a$ does not have the last property (to verify the latter claim, it suffices to take for $x_2$, $x_3$ and $x_4$ three pairwise incomparable elements of $L$). This example is communicated to the author by Shaprynski\v{\i}. 

Nevertheless, analogues of Lemma~\ref{join with neutr atom I-elem} and Corollary~\ref{join with SL I-elem} for cancellable elements are true. One can formulate these results explicitly.

\begin{lemma}[\!\!{\mdseries\cite[Lemma~2.1]{Gusev-Skokov-Vernikov-18}}]
\label{join with neutr atom canc}
Let $w$ be an atom and a neutral element of a lattice $L$. An element $x\in L$ is a cancellable element of $L$ if and only if the element $x\vee w$ has the same property.
\end{lemma}

This lemma and the properties of the variety $\mathbf{SL}$ mentioned above immediately imply the following

\begin{corollary}
\label{join with SL canc}
A \textup[commutative\textup] semigroup variety $\mathbf V$ is a cancellable element of the lattice $\mathbb{SEM}$ \textup[respectively $\mathbb{COM}$\textup] if and only if the variety $\mathbf{V\vee SL}$ has the same property.
\end{corollary}

\subsection{Modular and upper-modularl elements in lattices of equivalence relations}
\label{subsection Eq:mod,umod}

If $S$ is a set then $\Eq(S)$ stands for the lattice of equivalence relations on $S$.

\begin{proposition}
\label{Eq cmod,umod}
Let $S$ be a non-empty set. For an equivalence relation $\alpha$ on $S$, the following are equivalent:
\begin{itemize}
\item[\textup{a)}]$\alpha$ is a modular element of the lattice $\Eq(S)$;
\item[\textup{b)}]$\alpha$ is an upper-modular element of the lattice $\Eq(S)$;
\item[\textup{c)}]$\alpha$ has at most one non-singleton class.
\end{itemize}
\end{proposition}

The equivalences~a)\,$\Longleftrightarrow$\,c) and~b)\,$\Longleftrightarrow$\,c) of this proposition were proved in~\cite[Proposition~2.2]{Jezek-81} and~\cite[Proposition~3]{Vernikov-Volkov-88} respectively.

Proposition~\ref{Eq cmod,umod} turns out to be very helpful for the examination of modular and lower-modular elements in varietal lattices. In order to explain how this proposition can be applied, we need some new definitions and notation. Note that a semigroup $S$ satisfies the identity system $\mathbf wx\approx x\mathbf{w\approx w}$ where the letter $x$ does not occur in the word $\mathbf w$ if and only if $S$ contains a zero element~0 and all values of $\mathbf w$ in $S$ equal to~0. We adopt the usual convention of writing $\mathbf w\approx0$ as a short form of such a system and referring to the expression $\mathbf w\approx0$ as to a single identity. Identities of the form $\mathbf w\approx0$ are called 0-\emph{reduced}. Further, let $\mathbf X$ be a semigroup variety and $\mathbf{V\subseteq X}$, $F$ be the $\mathbf X$-free object and $\nu$ be the fully invariant congruence on $F$ corresponding to $\mathbf V$. It is clear that if $\mathbf V$ may be given within $\mathbf X$ by the family of 0-reduced identities $\{\mathbf w_i\approx0\mid i\in I\}$ only then $\nu$ has just one non-singleton class (namely, the collection of all equivalence classes modulo $\mathbf X$ that contain the words $\mathbf w_i$ where $i$ runs over $I$). Now recall the generally known fact that the lattice $L(\mathbf X)$ is anti-isomorphic to the lattice of all fully invariant congruences on $F$. Therefore, the lattice $\Eq(F)$ contains an anti-isomorphic copy of $L(\mathbf X)$. Finally, we note that the notion of a modular element is self-dual. Combining all these observations with Proposition~\ref{Eq cmod,umod}, we have the following

\begin{corollary}
\label{L(X) cmod,lmod suf}
Let $\mathbf X$ be a semigroup variety and $\mathbf{V\subseteq X}$. If $\mathbf V$ is defined within $\mathbf X$ by $0$-reduced identities only then $\mathbf V$ is a modular and lower-modular element of the lattice $L(\mathbf X)$.
\end{corollary}

\subsection{Special elements in subgroup lattices of finite symmetric groups}
\label{subsection Sub:all-but-umod,lmod}

The subgroup lattice of a group $G$ is denoted by $\Sub(G)$. We denote by $S_n$ the full symmetric group on the set $\{1,2,\dots,n\}$. If $\mathbf V$ is a semigroup variety and $n$ is a natural number then we put
$$
\Perm_n(\mathbf V)=\{\pi\in S_n\mid\mathbf V\text{ satisfies the identity }p_n[\pi]\}.
$$
Clearly, $\Perm_n(\mathbf V)$ is a subgroup in $S_n$. The following assertion explains our interest to modular and cancellable elements in lattices of the form $\Sub(S_n)$. 

\begin{lemma}
\label{Perm_n(V) cmod,canc}
Let $\mathbf V$ be a semigroup variety and $n$ be a natural number. If $\mathbf V$ is a modular \textup[cancellable\textup] element of the lattice $\mathbb{SEM}$ then the group $\Perm_n(\mathbf V)$ is a modular \textup[cancellable\textup] element of the lattice $\Sub(S_n)$.
\end{lemma}

This lemma is proved in~\cite[Corollary~4.3]{Vernikov-07-cmod} for modular elements and in~\cite[Proposition~3.2]{Shaprynskii-Skokov-Vernikov-19} for cancellable ones.

We denote by $A_n$ the alternative subgroup of $S_n$, by $V_4$ the Klein four-group and by $T$ the singleton group.

\begin{proposition}[\!\!{\mdseries\cite[Propositions~3.1,~3.7 and~3.8]{Jezek-81}}]
\label{Sub cmod}
Let $n$ be a natural number. A subgroup $G$ of a group $S_n$ is a modular element of the lattice $\Sub(S_n)$ if and only if one of the following holds:
\begin{itemize}
\item[\textup{(i)}]$n\le3$;
\item[\textup{(ii)}]$n=4$ and either $G=T$ or $G\supseteq V_4$;
\item[\textup{(iii)}]$n\ge5$ and $G$ coincides with either $T$ or $A_n$ or $S_n$.
\end{itemize}
\end{proposition}

This proposition is applied in the examination of modular elements in $\mathbb{SEM}$.

\begin{proposition}
\label{Sub neutr,stand,costand,canc,distr,codistr}
Let $n$ be a natural number. For a subgroup $G$ of the group $S_n$, the following are equivalent:
\begin{itemize}
\item[\textup{a)}]$G$ is a neutral element of the lattice $\Sub(S_n)$;
\item[\textup{b)}]$G$ is a standard element of the lattice $\Sub(S_n)$;
\item[\textup{c)}]$G$ is a costandard element of the lattice $\Sub(S_n)$;
\item[\textup{d)}]$G$ is a cancellable element of the lattice $\Sub(S_n)$;
\item[\textup{e)}]$G$ is a distributive element of the lattice $\Sub(S_n)$;
\item[\textup{f)}]$G$ is a codistributive element of the lattice $\Sub(S_n)$;
\item[\textup{g)}]either $G=T$ or $G=S_n$.
\end{itemize}
\end{proposition}

The equivalence of the claims~a)--c) and~e)--g) of this proposition is checked in the proof of~\cite[Theorem~2]{Vernikov-01}, while the equivalence~d)\,$\Longleftrightarrow$\,g) is verified in~\cite[Lemma~2.3]{Shaprynskii-Skokov-Vernikov-19}. Note that the key role in the proof of the equivalence~d)\,$\Longleftrightarrow$\,g) plays Proposition~\ref{Sub cmod}. Proposition~\ref{Sub neutr,stand,costand,canc,distr,codistr} is applied in the proof of Theorems~\ref{SEM canc} and~\ref{OC neutr,stand,costand,canc,distr,codistr}.

\subsection{Special elements in congruence lattices of $G$-sets}
\label{subsection gs:distr,codistr,canc,stand,costand,neutr}

A unary algebra with the carrier $A$ and the set of (unary) operations $G$ is called a $G$-\emph{set} if $G$ is equipped by a structure of a group and this group structure on $G$ is compatible with the unary structure on $A$ (this means that if $g,h\in G$, $x\in A$ and $e$ is the unit element of $G$ then $g\bigl(h(x)\bigr)=(gh)(x)$ and $e(x)=x$). Our interest to $G$-sets is explained by the fact that the lattice $\mathbb{OC}$ admits a concise and transparent description in terms of congruence lattices of $G$-sets. More precisely, $\mathbb{OC}$ is anti-isomorphic to a subdirect product of congruence lattices of countably infinite series of certain $G$-sets (see~\cite{Volkov-94} or~\cite[Subsection~5.1]{Shevrin-Vernikov-Volkov-09}). To apply this result for examination of special elements in $\mathbb{OC}$, some information about special elements in congruence lattices of $G$-sets is required.

A $G$-set $A$ is said to be \emph{transitive} if, for all $a,b\in A$, there exists $g\in G$ such that $g(a)=b$. If $A$ is a $G$-set and $a\in A$ then we put
$$
\Stab_A(a)=\bigl\{g\in G\mid g(a)=a\bigr\}.
$$
Clearly, $\Stab_A(a)$ is a subgroup in $G$. This subgroup is called a \emph{stabilizer} of an element $a$ in $A$. The congruence lattice of a $G$-set $A$ is denoted by $\Con(A)$.

\begin{proposition}
\label{G-sets neutr,stand,costand,canc,distr,codistr}
Let $A$ be a non-transitive $G$-set such that $\Stab_A(x)=\Stab_A(y)$ for all elements $x,y\in A$. For a congruence $\alpha$ on $A$, the following are equivalent:
\begin{itemize}
\item[\textup{a)}]$\alpha$ is a neutral element of the lattice $\Con(A)$;
\item[\textup{b)}]$\alpha$ is a standard element of the lattice $\Con(A)$;
\item[\textup{c)}]$\alpha$ is a costandard element of the lattice $\Con(A)$;
\item[\textup{d)}]$\alpha$ is a cancellable element of the lattice $\Con(A)$;
\item[\textup{e)}]$\alpha$ is a distributive element of the lattice $\Con(A)$;
\item[\textup{f)}]$\alpha$ is a codistributive element of the lattice $\Con(A)$;
\item[\textup{g)}]$\alpha$ is either the universal relation or the equality relation on $A$.
\end{itemize}
\end{proposition}

The equivalence of the claims~a)--c) and~e)--g) of this proposition was proved in~\cite[Theorem~1]{Vernikov-01}. The equivalence~d)\,$\Longleftrightarrow$\,~g) is verified in~\cite[Proposition~3.4]{Shaprynskii-Vernikov-21}. 

$G$-sets that appear in~\cite{Volkov-94} in the description of the lattice $\mathbb{OC}$ have the property that the stabilizer of any element in these $G$-sets is the trivial group. Thus, the application of Proposition~\ref{G-sets neutr,stand,costand,canc,distr,codistr} is not hindered by the hypothesis that stabilizers of all elements in $A$ coincide. It is presently unknown, whether the proposition holds without this hypothesis.

\subsection{Upper-modular and codistributive elements: interrelations between lattice identities and a hereditary property}
\label{subsection umod,codistr:ident&hered}

The following easy observation turns out to be helpful.

\begin{observation}
\label{abstr umod,codistr ident&hered}
Let $L$ be a lattice. If an element $a\in L$ is upper-modular \textup[codistributive\textup] in $L$ and the lattice $(a]$ is modular \textup[distributive\textup] then every element of $(a]$ is upper-modular \textup[codistributive\textup] in $L$.
\end{observation}

This claim was noted in~\cite[Lemma~2.1]{Vernikov-08-umod1} for upper-modular elements and in~\cite[Lemma~2.2]{Vernikov-11} for codistributive ones.

Observation~\ref{abstr umod,codistr ident&hered} immediately implies the following

\begin{corollary}
\label{SEM,COM umod,codistr ident&hered}
Let $\mathbf V$ be a \textup[commutative\textup] semigroup variety.
\begin{itemize}
\item[\textup{(i)}] If $\mathbf V$ is an upper-modular element of the lattice $\mathbb{SEM}$ \textup[respectively $\mathbb{COM}$\textup] and the lattice $L(\mathbf V)$ is modular then every subvariety of $\mathbf V$ is an upper-modular element of $\mathbb{SEM}$ \textup[respectively $\mathbb{COM}$\textup].
\item[\textup{(ii)}] If $\mathbf V$ is a codisributive element of the lattice $\mathbb{SEM}$ \textup[respectively $\mathbb{COM}$\textup] and the lattice $L(\mathbf V)$ is disributive then every subvariety of $\mathbf V$ is a codisributive element of $\mathbb{SEM}$ \textup[respectively $\mathbb{COM}$\textup].
\end{itemize}
\end{corollary}

\section{The lattice $\mathbb{SEM}$}
\label{section SEM}

For convenience, we call a semigroup variety \emph{modular} if it is a modular element of the lattice $\mathbb{SEM}$ and adopt similar agreement for all other types of special elements. The main results of this section provide:
\begin{itemize}
\item a complete classification of lower-modular, cancellable, distributive, standard, costandard or neutral varieties (Theorems~\ref{SEM lmod},~\ref{SEM distr,stand},~\ref{SEM costand,neutr} and~\ref{SEM canc}),{\sloppy

}
\item a classification of modular, upper-modular or codistributive varieties in some wide partial cases (Theorems~\ref{SEM cmod permut3},~\ref{SEM umod deg>2},~\ref{SEM umod stpermut} and~\ref{SEM codistr stpermut}),
\item strong necessary conditions for a semigroup variety to be modular, upper-modular or codistributive (Theorems~\ref{SEM cmod nec},~\ref{SEM cmod nil-nec},~\ref{SEM umod deg<3-nec} and~\ref{SEM codistr nec}),
\item a sufficient condition for a semigroup variety to be modular (Theorem~\ref{SEM cmod suf}).
\end{itemize}
One can mention also Proposition~\ref{I-varieties are periodic} that gives an important information about $I$-varieties and Proposition~\ref{SEM cmod permut} containing an interesting partial information about modular varieties.

\subsection{$I$-varieties}
\label{subsection SEM:I-varieties}

We denote by $\mathbf{SEM}$ the variety of all semigroups. A semigroup variety $\mathbf V$ is called \emph{proper} if $\mathbf{V\ne SEM}$.

The class of $I$-varieties includes all varieties with non-trivial identities in subvariety lattices (see Corollary~\ref{identities and I-varieties}). It follows from results of~\cite{Burris-Nelson-71-finite} that a semigroup variety $\mathbf V$ is periodic whenever the lattice $L(\mathbf V)$ satisfies some non-trivial identity. As we have already mentioned in Section~\ref{section intr}, results of the articles~\cite{Jezek-McKenzie-93,Vernikov-07-lmod,Vernikov-08-umod1} imply that if a proper variety $\mathbf V$ is a special element of one of the nine types we consider in $\mathbb{SEM}$ then it is periodic too. All these statements are generalized by the following

\begin{proposition}[\!\!{\mdseries\cite[Theorem~1]{Shaprynskii-12-I-elem}}]
\label{I-varieties are periodic}
A proper $I$-variety of semigroups is periodic.
\end{proposition}

Formally speaking, this proposition is not applicable to cancellable varieties because cancellable elements of a lattice are not $I$-elements. But, in actual fact, Proposition~\ref{I-varieties are periodic} implies that any proper cancellable variety is periodic because a cancellable variety is a modular one and modular elements are $I$-elements.

On the other hand, it is verified in~\cite[Theorem~2]{Shaprynskii-12-I-elem} that there are periodic varieties (moreover, nil-varieties) of semigroups that are not $I$-varieties. However, as far as we know, an explicit example of a periodic semigroup variety that is not an $I$-variety is unknown so far.

\subsection{Lower-modular varieties}
\label{subsection SEM:lmod}

Varieties that may be given by 0-reduced identities only are called 0-\emph{reduced}. We denote by $\mathbf T$ the trivial semigroup variety. A number of partial results concerning lower-modular varieties were obtained in~\cite{Vernikov-07-lmod,Vernikov-08-lmod,Vernikov-Volkov-88}. All of them are covered by the following

\begin{theorem}
\label{SEM lmod}
A semigroup variety $\mathbf V$ is lower-modular if and only if either $\mathbf{V=SEM}$ or $\mathbf{V=M\vee N}$ where $\mathbf M$ is one of the varieties $\mathbf T$ or $\mathbf{SL}$, while $\mathbf N$ is a $0$-reduced variety.
\end{theorem}

This theorem was verified for the first time in~\cite[Theorem~1.1]{Shaprynskii-Vernikov-10} and was reproved in a simpler and shorter way in~\cite{Shaprynskii-12-mod-lmod}. The proof of Theorem~\ref{SEM lmod} given in~\cite{Shaprynskii-12-mod-lmod} is based on Theorem~\ref{L(oc) cmod,lmod per-nec} below. Note that the `if' part of Theorem~\ref{SEM lmod} immediately follows from Corollaries~\ref{L(X) cmod,lmod suf} (with $\mathbf{X=SEM}$) and~\ref{join with SL I-elem}.

Neutral, standard and distributive varieties are lower-modular. In view of Theorem~\ref{SEM lmod}, a description of varieties of these three types should look as follows:

\emph{A semigroup variety $\mathbf V$ is distributive \textup[standard, neutral\textup] if and only if either $\mathbf{V=SEM}$ or $\mathbf{V=M\vee N}$ where $\mathbf M$ is one of the varieties $\mathbf T$ or $\mathbf{SL}$, while $\mathbf N$ is a $0$-reduced variety such that \dots} (with some additional restriction to $\mathbf N$ depending on the type of element we consider).

Exact formulations of corresponding results are given in the following two subsections.

\subsection{Distributive and standard varieties}
\label{subsection SEM:distr,stand}

Put
\begin{align*}
&\mathbf Q=\var\{x^2y\approx xyx\approx yx^2\approx0\},\\
&\mathbf Q_n=\var\{x^2y\approx xyx\approx yx^2\approx x_1x_2\cdots x_n\approx0\},\\ 
&\mathbf R=\var\{x^2\approx xyx\approx0\},\\
&\mathbf R_n=\var\{x^2\approx xyx\approx x_1x_2\cdots x_n\approx0\}
\end{align*}
where $n\ge2$. It is easy to see that varieties of these four types are precisely all non-trivial 0-reduced varieties satisfying the identities $ x^2y\approx xyx\approx yx^2\approx0$.

\begin{theorem}
\label{SEM distr,stand}
For a semigroup variety $\mathbf V$, the following are equivalent:
\begin{itemize}
\item[\textup{a)}]$\mathbf V$ is distributive;
\item[\textup{b)}]$\mathbf V$ is standard;
\item[\textup{c)}]either $\mathbf{V=SEM}$ or $\mathbf{V=M\vee N}$ where $\mathbf M$ is one of the varieties $\mathbf T$ or $\mathbf{SL}$, while $\mathbf N$ is one of the varieties $\mathbf T$, $\mathbf Q$, $\mathbf Q_n$, $\mathbf R$ or $\mathbf R_n$.
\end{itemize}
\end{theorem}

The equivalence~a)\,$\Longleftrightarrow$\,c) is proved in~\cite[Theorem~1.1]{Vernikov-Shaprynskii-10}. Note that the proof of the implication~a)\,$\Longrightarrow$\,c) given in~\cite{Vernikov-Shaprynskii-10} may be essentially simplified by using Theorem~\ref{SEM lmod}. The implication~b)\,$\Longrightarrow$\,a) is evident. To verify the implication~a)\,$\Longrightarrow$\,b), we need two ingredients. First, it is verified in~\cite[Corollary~1.2]{Vernikov-Shaprynskii-10} that a distributive semigroup variety is a modular one\footnote{In afctually fact, this claim immediately follows from the implication~a)\,$\Longrightarrow$\,c) of Theorem~\ref{SEM distr,stand} and Corollaries~\ref{L(X) cmod,lmod suf} (with $\mathbf{X=SEM}$) and~\ref{join with SL I-elem}.}. Second, it is well known that an element of a lattice is standard whenever it is simultaneously distributive and modular (see~\cite[Lemma~II.1.1]{Gratzer-Schmidt-61}, for instance). Note that the former statement is strengthened by Corollary~\ref{SEM lmod is cmod} below.

It is well known that the set of all standard elements of an arbitrary lattice $L$ forms a sublattice in $L$~\cite[Theorem~259 on p.\,230]{Gratzer-11}. Theorem~\ref{SEM distr,stand} shows that the lattice of all distributive (equivalently, standard) varieties has the form shown in Fig.~\ref{SEM distr,stand diagram}.

\begin{figure}[tbh]
\begin{center}
\unitlength=1mm
\linethickness{0.4pt}
\begin{picture}(61,100)
\put(16,35){\circle*{1.33}}
\put(16,45){\circle*{1.33}}
\put(16,55){\circle*{1.33}}
\put(16,75){\circle*{1.33}}
\put(26,5){\circle*{1.33}}
\put(26,15){\circle*{1.33}}
\put(26,25){\circle*{1.33}}
\put(26,35){\circle*{1.33}}
\put(26,45){\circle*{1.33}}
\put(26,65){\circle*{1.33}}
\put(36,45){\circle*{1.33}}
\put(36,55){\circle*{1.33}}
\put(36,65){\circle*{1.33}}
\put(36,85){\circle*{1.33}}
\put(46,15){\circle*{1.33}}
\put(46,25){\circle*{1.33}}
\put(46,35){\circle*{1.33}}
\put(46,45){\circle*{1.33}}
\put(46,55){\circle*{1.33}}
\put(46,75){\circle*{1.33}}
\put(56,95){\circle*{1.33}}
\gasset{AHnb=0,linewidth=0.4}
\drawline[dash={0.3 1.5}{0}](16,55)(17,75)
\drawline[dash={0.3 1.5}{0}](26,45)(26,65)
\drawline[dash={0.3 1.5}{0}](36,65)(36,85)
\drawline[dash={0.3 1.5}{0}](46,55)(46,75)
\drawline(26,25)(26,35)(16,45)(36,55)(46,45)(26,35)
\drawline(56,95)(16,75)(26,65)(46,75)(36,85)
\drawline(26,25)(46,35)(36,45)(16,35)
\drawline(26,35)(26,45)(46,55)
\drawline(16,55)(26,45)
\drawline(26,15)(46,25)
\drawline(36,45)(36,65)
\drawpolygon(26,5)(26,25)(16,35)(16,55)(36,65)(46,55)(46,15)
\small
\put(16,78){\makebox(0,0)[cc]{$\mathbf Q$}}
\put(25,15){\makebox(0,0)[rc]{$\mathbf Q_2=\mathbf R_2=\mathbf{ZM}$}}
\put(15,35){\makebox(0,0)[rc]{$\mathbf Q_3$}}
\put(15,45){\makebox(0,0)[rc]{$\mathbf Q_4$}}
\put(15,55){\makebox(0,0)[rc]{$\mathbf Q_5$}}
\put(26,68){\makebox(0,0)[cc]{$\mathbf R$}}
\put(27,24){\makebox(0,0)[lc]{$\mathbf R_3$}}
\put(27,34){\makebox(0,0)[lc]{$\mathbf R_4$}}
\put(27,44){\makebox(0,0)[lc]{$\mathbf R_5$}}
\put(47,15){\makebox(0,0)[lc]{$\mathbf{SL}$}}
\put(56,98){\makebox(0,0)[cc]{$\mathbf{SEM}$}}
\put(26,2){\makebox(0,0)[cc]{$\mathbf T$}}
\end{picture}
\caption{The lattice of distributive varieties}
\label{SEM distr,stand diagram}
\end{center}
\end{figure}
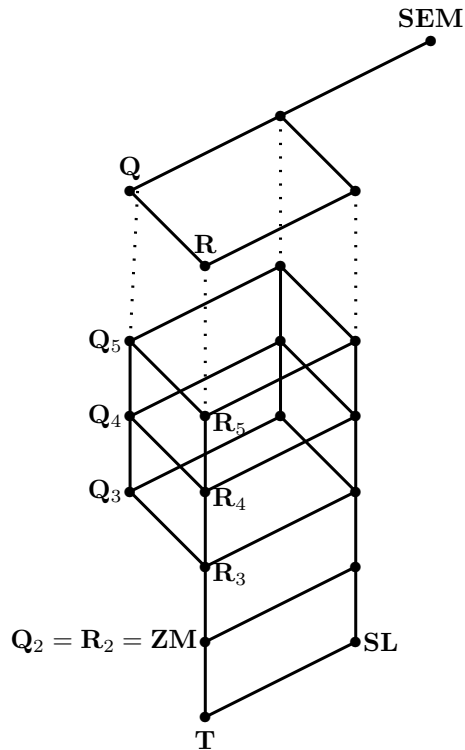

Thus, we have

\begin{corollary}
\label{SEM distr,stand sublat}
The class of all distributive \textup[standard\textup] semigroup varieties forms a countably infinite distributive sublattice of the lattice $\mathbb{SEM}$.
\end{corollary}

\subsection{Costandard and neutral varieties}
\label{subsection SEM:costand,neutr}

Results of the earlier works~\cite{Aizenstat-74,Melnik-69,Melnik-70,Melnik-73,Salij-69} mentioned in the Introduction easily imply that the varieties $\mathbf{SL}$ and $\mathbf{ZM}$ are neutral. Since the set of all neutral elements of a lattice $L$ forms a sublattice of $L$ (see~\cite[Theorem~259(iii) on p.\,230]{Gratzer-11}), the variety $\mathbf{SL\vee ZM}$ is neutral too. The neutrality of the varieties $\mathbf T$ and $\mathbf{SEM}$ is evident. It turns out that there are no neutral varieties except the five mentioned ones. More exactly, the following statement is true.

\begin{theorem}
\label{SEM costand,neutr}
For a semigroup variety $\mathbf V$, the following are equivalent:
\begin{itemize}
\item[\textup{a)}]$\mathbf V$ is both lower-modular and upper-modular;
\item[\textup{b)}]$\mathbf V$ is both distributive and codistributive;
\item[\textup{c)}]$\mathbf V$ is costandard;
\item[\textup{d)}]$\mathbf V$ is neutral;
\item[\textup{e)}]$\mathbf V$ is one of the varieties $\mathbf T$, $\mathbf{SL}$, $\mathbf{ZM}$, $\mathbf{SL\vee ZM}$ or $\mathbf{SEM}$.
\end{itemize}
\end{theorem}

Clearly, the claim~e) of this theorem may be reformulated in the manner specified in Subsection~\ref{subsection SEM:lmod}: either $\mathbf{V=SEM}$ or $\mathbf{V=M\vee N}$ where $\mathbf M$ is one of the varieties $\mathbf T$ or $\mathbf{SL}$, while $\mathbf N$ is a 0-reduced variety such that $xy\approx0$ in $\mathbf N$.

We do not include in Theorem~\ref{SEM costand,neutr} the claim that $\mathbf V$ is both standard and costandard because it is well known that an element of arbitrary lattice is both standard and costandard if and only if it is neutral (see~\cite[Theorem~255 on p.\,228]{Gratzer-11}, for instance). The equivalences a)\,$\Longleftrightarrow$\,~e),~c)\,$\Longleftrightarrow$\,e) and d)\,$\Longleftrightarrow$\,e) of Theorem~\ref{SEM costand,neutr} were verified in~\cite[Corollary~3.5]{Vernikov-08-lmod},~\cite[Theorem~1.3]{Vernikov-11} and~\cite[Proposition~4.1]{Volkov-05} respectively, while the implications~d)\,$\Longrightarrow$\,b)\,$\Longrightarrow$\,a) are evident.

Since a neutral element of a lattice is standard, Theorem~\ref{SEM costand,neutr} implies the following

\begin{corollary}[\!\!{\mdseries\cite[Corollary~1.1]{Vernikov-11}}]
\label{SEM costand is stand}
Every costandard semigroup variety is standard.
\end{corollary}

\subsection{An application to definable varieties}
\label{subsection SEM:definability}

Here we discuss an interesting application of results overviewed above. We need some new definitions. A subset $A$ of a lattice $\langle L;\vee,\wedge\rangle$ is called \emph{definable in} $L$ if there exists a first-order formula $\Phi(x)$ with one free variable $x$ in the language of lattice operations $\vee$ and $\wedge$ which \emph{defines $A$ in} $L$. This means that, for an element $a\in L$, the sentence $\Phi(a)$ is true if and only if $a\in A$. If $A$ consists of a single element, then we talk about definability of this element. A set $\mathcal X$ of semigroup varieties (or a single semigroup variety $\mathbf X$) is said to be \emph{definable} if it is definable in $\mathbb{SEM}$. In this situation we will say that the corresponding first-order formula \emph{defines} the set $\mathcal X$ or the variety $\mathbf X$.

A number of deep results about definable varieties and sets of varieties of semigroups have been obtained in~\cite{Jezek-McKenzie-93} by Je\v{z}ek and McKenzie. It has been conjectured there that every finitely based semigroup variety is definable up to duality. The conjecture is confirmed in~\cite{Jezek-McKenzie-93} for locally finite finitely based varieties. On their way to obtaining this fundamental result, Je\v{z}ek and McKenzie proved the definability of several important sets of semigroup varieties such as the sets of all finitely based, all locally finite, all finitely generated and all 0-reduced semigroup varieties. But the article~\cite{Jezek-McKenzie-93} contains no explicit first-order formulas that define any of these sets of varieties. The task of writing an explicit formula that defines the set of all finitely based or the set of all locally finite or the set of all finitely generated varieties seems to be extremely difficult. On the other hand, the set of all 0-reduced varieties can be defined by a quite simple first-order formula based on descriptions of lower-modular and neutral varieties.

Indeed, Theorem~\ref{SEM lmod} shows that a semigroup variety is 0-reduced if and only if it is lower-modular and does not contain the variety $\mathbf{SL}$. It remains to define the variety $\mathbf{SL}$. Theorem~\ref{SEM costand,neutr} together with the well-known description of atoms of the lattice $\mathbb{SEM}$ (see~\cite[Section~1]{Shevrin-Vernikov-Volkov-09}, for instance) imply that this lattice contains exactly two neutral atoms, namely the varieties $\mathbf{SL}$ and $\mathbf{ZM}$. Recall that a semigroup variety $\mathbf V$ is called \emph{chain} if the lattice $L(\mathbf V)$ is a chain. It is well known that the variety $\mathbf{ZM}$ is properly contained in some chain variety, while the variety $\mathbf{SL}$ is not (see~\cite{Sukhanov-82} or~\cite{Vernikov-12} for more details). Combining the mentioned observations, we see that the set of all 0-reduced varieties may be defined as the set $\mathcal K$ of semigroup varieties with the following properties:
\begin{itemize}
\item[(i)]every member of $\mathcal K$ is a lower-modular variety;
\item[(ii)]if $\mathbf V\in\mathcal K$ and $\mathbf V$ contains some neutral atom $\mathbf A$ then $\mathbf A$ is properly contained in some chain variety.
\end{itemize}
It is evident that properties~(i) and~(ii) may be written by simple first-order formulas with one free variable.

An explicit formula that defines the set of all 0-reduced varieties is written in~\cite[Section~3]{Vernikov-12}. Note that the description of distributive semigroup varieties given by Theorem~\ref{SEM distr,stand} may also be applied to define some interesting varieties (see~\cite[Section~6]{Vernikov-12}).

\subsection{Modular varieties}
\label{subsection SEM:cmod}

The problem of description of modular semigroup varieties is open so far. Here we provide some partial results concerning this problem.

Recall that a semigroup variety is called a \emph{nil-variety} if it consists of nilsemigroups or, equivalently, satisfies an identity of the form $x^n\approx0$ for some natural $n$. Clearly, every 0-reduced variety is a nil-variety. The following theorem gives a strong necessary condition for a semigroup variety to be modular.

\begin{theorem}
\label{SEM cmod nec}
If $\mathbf V$ is a modular semigroup variety then either $\mathbf{V=SEM}$ or $\mathbf{V=M\vee N}$ where $\mathbf M$ is one of the varieties $\mathbf T$ or $\mathbf{SL}$, while $\mathbf N$ is a nil-variety.
\end{theorem}

This theorem readily follows from~\cite[Proposition~1.6]{Jezek-McKenzie-93}. A deduction of Theorem~\ref{SEM cmod nec} from~\cite[Proposition~1.6]{Jezek-McKenzie-93} is given explicitly in~\cite[Proposition~2.1]{Vernikov-07-cmod}. A direct and transparent proof of Theorem~\ref{SEM cmod nec} not depending on a technique from~\cite{Jezek-McKenzie-93} is given in~\cite{Shaprynskii-12-mod-lmod}. This proof is based on Theorem~\ref{L(oc) cmod,lmod per-nec} below.

Theorem~\ref{SEM cmod nec} and Corollary~\ref{join with SL I-elem} completely reduce the examination of modular varieties to nil-varieties. There is a strong necessary condition for a nil-variety to be modular. To formulate this result, we need some additional definitions.

We call an identity $\mathbf{u\approx v}$ \emph{substitutive} if the words $\mathbf u$ and $\mathbf v$ depend on the same letters and $\mathbf v$ may be obtained from $\mathbf u$ by renaming of letters. In~\cite{Jezek-81}, Je\v{z}ek describes modular elements of the lattice of all varieties (more precisely, all equational theories) of any given type. In particular, it follows from~\cite[Lemma~6.3]{Jezek-81} that if a nil-variety of semigroups $\mathbf V$ is a modular element of the lattice of all groupoid varieties then $\mathbf V$ may be given by 0-reduced and substitutive identities only. This does not imply directly the same conclusion for modular nil-varieties because a modular element of $\mathbb{SEM}$ need not be a modular element of the lattice of all groupoid varieties. Nevertheless, the following assertion shows that the `semigroup analogue' of the mentioned result of Je\v{z}ek holds true.

\begin{theorem}[\!\!{\mdseries\cite[Proposition~2.2]{Vernikov-07-cmod}}]
\label{SEM cmod nil-nec}
A modular nil-variety of semigroups may be given by $0$-reduced and substitutive identities only.
\end{theorem}

Corollary~\ref{L(X) cmod,lmod suf} with $\mathbf{X=SEM}$ immediately implies the following

\begin{theorem}
\label{SEM cmod suf}
Every $0$-reduced semigroup variety is modular.
\end{theorem}

This fact was noted for the first time in~\cite[Corollary~3]{Vernikov-Volkov-88} and rediscovered (in other terminology) in~\cite[Proposition~1.1]{Jezek-McKenzie-93}\footnote{In fact, Theorem~\ref{SEM cmod suf} is an immediate consequence of the equivalence of the claims~a) and~c) of Proposition~\ref{Eq cmod,umod}. This equivalence is proved by Je\v{z}ek in~\cite[Proposition~2.2]{Jezek-81}. In~\cite{Vernikov-Volkov-88}, Theorem~\ref{SEM cmod suf} is justified just by reference to this result by Je\v{z}ek. It is rather unexpected that Je\v{z}ek and McKenzie give much more complex proof of Theorem~\ref{SEM cmod suf} in~\cite{Jezek-McKenzie-93}.}.

Theorems~\ref{SEM cmod nil-nec} and~\ref{SEM cmod suf} provide a necessary and a sufficient condition for a nil-variety to be modular respectively. The gap between these conditions seems to be not very large. But the necessary condition is not a sufficient one, while the sufficient condition is not a necessary one (this follows from Corollary~\ref{SEM cmod,canc commut} below).

Theorems~\ref{SEM lmod} and~\ref{SEM cmod suf} and Corollary~\ref{join with SL I-elem} immediately imply the following

\begin{corollary}[\!\!{\mdseries\cite[Corollary~1.2]{Shaprynskii-Vernikov-10}}]
\label{SEM lmod is cmod}
Every lower-modular semigroup variety is modular.
\end{corollary}

Theorems~\ref{SEM cmod nil-nec} and~\ref{SEM cmod suf} show that in order to describe modular nil-varieties (and therefore, all modular varieties) we need to examine nil-varieties satisfying substitutive identities. A natural partial case of substitutive identities are permutational ones. Some interesting information concerning modular varieties satisfying a permutational identity give the following assertion.

\begin{proposition}[\!\!{\mdseries\cite[Theorem~4.5]{Vernikov-07-cmod}}]
\label{SEM cmod permut}
Let $\mathbf V$ be a modular variety of semigroups satisfying a permutational identity $p_n[\pi]$. Then $\mathbf V$ satisfies also:
\begin{itemize}
\item[\textup{(i)}]all permutational identities of length $n+1$;
\item[\textup{(ii)}]all permutational identities of length $n$ whenever $n\ge5$ and the permutation $\pi$ is odd;
\item[\textup{(iii)}]an arbitrary identity of the form $\mathbf u\approx0$, where $\mathbf u$ is a word of length $n$ depending on $n-1$ letters, whenever $n\ge4$ and $\mathbf V$ is a nil-variety.
\end{itemize}
\end{proposition}

The proof of this proposition is based on Theorem~\ref{SEM cmod nil-nec}, Proposition~\ref{Sub cmod}, Lemma~\ref{Perm_n(V) cmod,canc} and results of the article~\cite{Pollak-73}. 

The following statement gives a complete classification of modular varieties satisfying a permutational identity of length~3. 

\begin{theorem}[\!\!{\mdseries\cite[Theorem~1.1]{Skokov-Vernikov-19}}]
\label{SEM cmod permut3}
A semigroup variety $\mathbf V$ satisfying a permutational identity of length~$3$ is modular if and only if $\mathbf{V=M\vee N}$ where $\mathbf M$ is one of the varieties $\mathbf T$ or $\mathbf{SL}$, while the variety $\mathbf N$ satisfies one of the following identity systems:
\begin{align*}
&xyz\approx zyx,\,x^2y\approx0;\\
&xyz\approx yzx,\,x^2y\approx0;\\
&xyz\approx yxz,\,xyzt\approx xzty,\,xy^2\approx0;\\
&xyz\approx xzy,\,xyzt\approx yzxt,\,x^2y\approx0.
\end{align*}
\end{theorem}

\subsection{Cancellable varieties}
\label{subsection SEM:canc}

Some partial results about cancellable semigroup varieties are proved in the articles~\cite{Gusev-Skokov-Vernikov-18,Skokov-Vernikov-19}. All these results are covered by the complete description of cancellable varieties obtained in~\cite[Theorem~1.1]{Shaprynskii-Skokov-Vernikov-19}. To formulate this result, we need notation for varieties introduced in Subsection~\ref{subsection SEM:distr,stand} and also some new notation. If $\mathbf X$ is one of the varieties $\mathbf Q$, $\mathbf Q_n$, $\mathbf R$ or $\mathbf R_n$, and $m\ge 2$ then we denote by $\mathbf X[m]$ the subvariety of $\mathbf X$ given within $\mathbf X$ by all permutational identities of length $m$. In particular, $\mathbf Q_n[m]=\mathbf Q_n$ and $\mathbf R_n[m]=\mathbf R_n$ whenever $m\ge n$.

\begin{theorem}[\!\!{\mdseries\cite[Theorem~1.1]{Shaprynskii-Skokov-Vernikov-19}}]
\label{SEM canc}
A semigroup variety $\mathbf V$ is cancellable if and only if either $\mathbf{V=SEM}$ or $\mathbf{V=M\vee N}$ where $\mathbf M$ is one of the varieties $\mathbf T$ or $\mathbf{SL}$, while $\mathbf N$ is one of the varieties $\mathbf T$, $\mathbf Q$, $\mathbf R$, $\mathbf Q[m]$ with $m\ge 2$, $\mathbf R[m]$ with $m\ge 2$, $\mathbf Q_n[m]$ with $2\le m\le n$ or $\mathbf R_n[m]$ with $2\le m\le n$.
\end{theorem}

This theorem shows that if a cancellable semigroup variety satisfies a permutational identity of length $m$ then it satisfies all permutational identities of length $m$ (see~\cite[Proposition~3.2]{Shaprynskii-Skokov-Vernikov-19}). In fact, this claim follows from Lemma~\ref{Perm_n(V) cmod,canc} and Proposition~\ref{Sub neutr,stand,costand,canc,distr,codistr}. Theorem~\ref{SEM canc} shows also that cancellable nil-varieties satisfies the following claim which stronger than Theorem~\ref{SEM cmod nil-nec}: such varieties can be given by 0-reduced and permutational identities only.

It is known that the set of all cancellable elements of a lattice $L$ may not be a sublattice of $L$. However, Theorem~\ref{SEM distr,stand} readily implies the following analogue of Corollary~\ref{SEM distr,stand sublat}.

\begin{corollary}
\label{SEM canc sublat}
The class of all cancellable semigroup varieties forms a countably infinite distributive sublattice of the lattice $\mathbb{SEM}$.
\end{corollary}

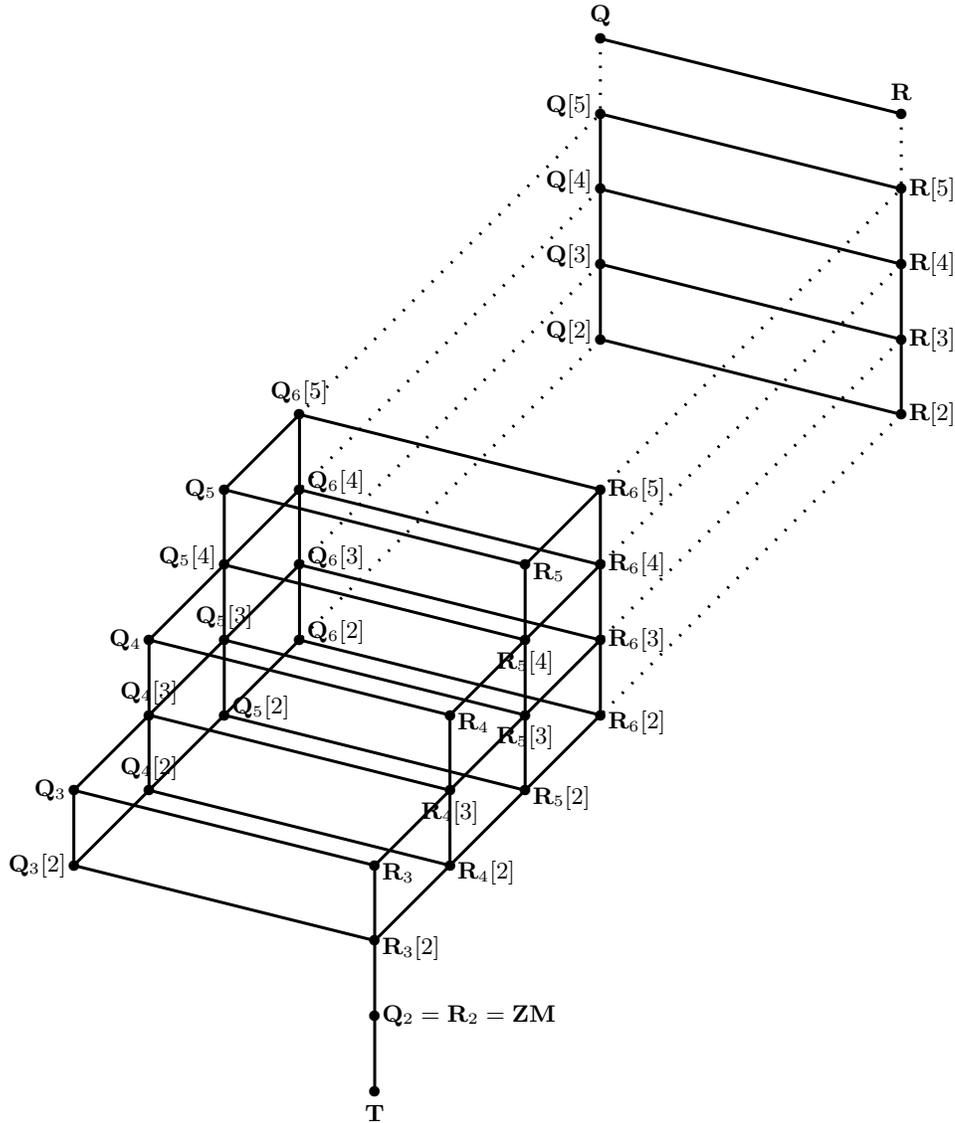
\begin{figure}[tbh]
\begin{center}
\unitlength=1mm
\linethickness{0.4pt}
\begin{picture}(127,150)
\put(9,35){\circle*{1.33}}
\put(9,45){\circle*{1.33}}
\put(19,45){\circle*{1.33}}
\put(19,55){\circle*{1.33}}
\put(19,65){\circle*{1.33}}
\put(29,55){\circle*{1.33}}
\put(29,65){\circle*{1.33}}
\put(29,75){\circle*{1.33}}
\put(29,85){\circle*{1.33}}
\put(39,65){\circle*{1.33}}
\put(39,75){\circle*{1.33}}
\put(39,85){\circle*{1.33}}
\put(39,95){\circle*{1.33}}
\put(49,5){\circle*{1.33}}
\put(49,15){\circle*{1.33}}
\put(49,25){\circle*{1.33}}
\put(49,35){\circle*{1.33}}
\put(59,35){\circle*{1.33}}
\put(59,45){\circle*{1.33}}
\put(59,55){\circle*{1.33}}
\put(69,45){\circle*{1.33}}
\put(69,55){\circle*{1.33}}
\put(69,65){\circle*{1.33}}
\put(69,75){\circle*{1.33}}
\put(79,55){\circle*{1.33}}
\put(79,65){\circle*{1.33}}
\put(79,75){\circle*{1.33}}
\put(79,85){\circle*{1.33}}
\put(79,105){\circle*{1.33}}
\put(79,115){\circle*{1.33}}
\put(79,125){\circle*{1.33}}
\put(79,135){\circle*{1.33}}
\put(79,145){\circle*{1.33}}
\put(119,95){\circle*{1.33}}
\put(119,105){\circle*{1.33}}
\put(119,115){\circle*{1.33}}
\put(119,125){\circle*{1.33}}
\put(119,135){\circle*{1.33}}
\gasset{AHnb=0,linewidth=0.4}
\drawline[dash={0.3 1.5}{0}](39,65)(79,105)
\drawline[dash={0.3 1.5}{0}](39,75)(79,115)
\drawline[dash={0.3 1.5}{0}](39,85)(79,125)
\drawline[dash={0.3 1.5}{0}](39,95)(79,135)(79,145)
\drawline[dash={0.3 1.5}{0}](79,55)(119,95)
\drawline[dash={0.3 1.5}{0}](79,65)(119,105)
\drawline[dash={0.3 1.5}{0}](79,75)(119,115)
\drawline[dash={0.3 1.5}{0}](79,85)(119,125)(119,135)
\drawline(49,5)(49,35)(9,45)(9,35)(49,25)(79,55)(79,85)(39,95)(39,65)(9,35)
\drawline(9,45)(39,75)(79,65)(49,35)
\drawline(29,55)(69,45)(69,75)(29,85)
\drawline(29,55)(29,85)(39,95)
\drawline(19,55)(59,45)
\drawline(19,65)(59,55)
\drawline(29,65)(69,55)
\drawline(29,75)(69,65)
\drawline(39,65)(79,55)
\drawline(69,75)(79,85)
\drawline(79,115)(119,105)
\drawline(79,125)(119,115)
\drawline(79,145)(119,135)
\drawpolygon(19,45)(19,65)(39,85)(79,75)(59,55)(59,35)
\drawpolygon(79,105)(79,135)(119,125)(119,95)
\footnotesize
\put(79,148){\makebox(0,0)[cc]{$\mathbf Q$}}
\put(50,15){\makebox(0,0)[lc]{$\mathbf Q_2=\mathbf R_2=\mathbf{ZM}$}}
\put(78,106){\makebox(0,0)[rc]{$\mathbf Q[2]$}}
\put(78,116){\makebox(0,0)[rc]{$\mathbf Q[3]$}}
\put(8,45){\makebox(0,0)[rc]{$\mathbf Q_3$}}
\put(8,35){\makebox(0,0)[rc]{$\mathbf Q_3[2]$}}
\put(78,126){\makebox(0,0)[rc]{$\mathbf Q[4]$}}
\put(18,65){\makebox(0,0)[rc]{$\mathbf Q_4$}}
\put(19,48){\makebox(0,0)[cc]{$\mathbf Q_4[2]$}}
\put(19,58){\makebox(0,0)[cc]{$\mathbf Q_4[3]$}}
\put(78,136){\makebox(0,0)[rc]{$\mathbf Q[5]$}}
\put(28,85){\makebox(0,0)[rc]{$\mathbf Q_5$}}
\put(30,56){\makebox(0,0)[lc]{$\mathbf Q_5[2]$}}
\put(29,68){\makebox(0,0)[cc]{$\mathbf Q_5[3]$}}
\put(28,76){\makebox(0,0)[rc]{$\mathbf Q_5[4]$}}
\put(40,66){\makebox(0,0)[lc]{$\mathbf Q_6[2]$}}
\put(40,76){\makebox(0,0)[lc]{$\mathbf Q_6[3]$}}
\put(40,86){\makebox(0,0)[lc]{$\mathbf Q_6[4]$}}
\put(39,98){\makebox(0,0)[cc]{$\mathbf Q_6[5]$}}
\put(119,138){\makebox(0,0)[cc]{$\mathbf R$}}
\put(120,95){\makebox(0,0)[lc]{$\mathbf R[2]$}}
\put(120,105){\makebox(0,0)[lc]{$\mathbf R[3]$}}
\put(50,34){\makebox(0,0)[lc]{$\mathbf R_3$}}
\put(50,24){\makebox(0,0)[lc]{$\mathbf R_3[2]$}}
\put(120,115){\makebox(0,0)[lc]{$\mathbf R[4]$}}
\put(60,54){\makebox(0,0)[lc]{$\mathbf R_4$}}
\put(60,34){\makebox(0,0)[lc]{$\mathbf R_4[2]$}}
\put(59,42){\makebox(0,0)[cc]{$\mathbf R_4[3]$}}
\put(120,125){\makebox(0,0)[lc]{$\mathbf R[5]$}}
\put(70,74){\makebox(0,0)[lc]{$\mathbf R_5$}}
\put(70,44){\makebox(0,0)[lc]{$\mathbf R_5[2]$}}
\put(69,52){\makebox(0,0)[cc]{$\mathbf R_5[3]$}}
\put(69,62){\makebox(0,0)[cc]{$\mathbf R_5[4]$}}
\put(80,54){\makebox(0,0)[lc]{$\mathbf R_6[2]$}}
\put(80,65){\makebox(0,0)[lc]{$\mathbf R_6[3]$}}
\put(80,75){\makebox(0,0)[lc]{$\mathbf R_6[4]$}}
\put(80,85){\makebox(0,0)[lc]{$\mathbf R_6[5]$}}
\put(49,2){\makebox(0,0)[cc]{$\mathbf T$}}
\end{picture}
\caption{The lattice of cancellable nil-varieties}
\label{SEM canc diagram}
\end{center}
\end{figure}

In Fig.~\ref{SEM canc diagram} we show the `main part' of this lattice, namely, the lattice of cancellable nil-varieties. To get the whole lattice of cancellable varieties, one need to adjoin the new greatest element (the variety $\mathbf{SEM}$) to the direct product of the lattice shown in Fig.~\ref{SEM canc diagram} and the 2-element chain (consisting of the varieties $\mathbf T$ and $\mathbf{SL}$).

Theorems~\ref{SEM cmod permut3} and~\ref{SEM canc} readily imply the following

\begin{corollary}
\label{SEM cmod,canc commut}
For a commutative semigroup variety $\mathbf V$, the following are equivalent:
\begin{itemize}
\item[\textup{a)}]$\mathbf V$ is modular;
\item[\textup{b)}]$\mathbf V$ is cancellable;
\item[\textup{c)}]$\mathbf{V=M\vee N}$ where $\mathbf M$ is one of the varieties $\mathbf T$ or $\mathbf{SL}$, while $\mathbf N$ is a commutative variety satisfying the identity
\begin{equation}
\label{xxy=0}
x^2y\approx0.
\end{equation}
\end{itemize}
\end{corollary}

The equivalences a)\,$\Longleftrightarrow$\,c) and b)\,$\Longleftrightarrow$\,c) in this corollary were first proved in~\cite[Theorem~3.1]{Vernikov-07-cmod} and~\cite[Theorem~1.1]{Gusev-Skokov-Vernikov-18} respectively. Note that the part of Corollary~\ref{SEM cmod,canc commut} concerning modular varieties is reproduced in~\cite{Vernikov-15} with inaccuracy\footnote{Namely, it is written in~\cite[Theorem~3.10]{Vernikov-15} that a commutative semigroup variety is modular if and only if it satisfies the identity~\eqref{xxy=0}.}.

\subsection{Upper-modular varieties}
\label{subsection SEM:umod}

The problem of description of upper-modular semigroup varieties is open so far. Here we provide some partial results concerning this problem. The first result classifies upper-modular varieties in some wide class of varieties. To formulate this statement we need some additional definitions and notation.

A semigroup variety $\mathbf V$ is called a variety of \emph{finite degree} [a variety of \emph{degree} $n$] if all nilsemigroups in $\mathbf V$ are nilpotent [if nilpotency degrees of nilsemigroups in $\mathbf V$ are bounded by the number $n$ and $n$ is the least number with this property]. We say that a semigroup variety is a \emph{variety of degree} $>n$ if it is either a variety of a finite degree $m$ with $m>n$ or not a variety of finite degree. Put
\begin{align*}
&\mathbf A_n=\var\{x^ny\approx y,\,xy\approx yx\}\enskip\text{where}\enskip n\ge1,\\
&\mathbf C=\var\{x^2\approx x^3,\,xy\approx yx\}.
\end{align*}
In particular, $\mathbf A_1=\mathbf T$. Note that $\mathbf A_n$ is the variety of all Abelian groups of exponent $n$.

\begin{theorem}[\!\!{\mdseries\cite[Theorem~1]{Vernikov-08-umod2}}]
\label{SEM umod deg>2}
A semigroup variety $\mathbf V$ of degree $>2$ is upper-modular if and only if one of the following holds:
\begin{itemize}
\item[\textup{(i)}]$\mathbf{V=SEM}$;
\item[\textup{(ii)}]$\mathbf{V=M\vee N}$ where $\mathbf M$ is one of the varieties $\mathbf T$ or $\mathbf{SL}$, while $\mathbf N$ is a nil-variety satisfying the identities
\begin{equation}
\label{xxy=xyy,xy=yx}
x^2y\approx xy^2,\,xy\approx yx;
\end{equation}
\item[\textup{(iii)}]$\mathbf{V=A}_n\vee\mathbf{M\vee N}$ where $n\ge1$, $\mathbf M$ is one of the varieties $\mathbf T$, $\mathbf{SL}$ or $\mathbf C$, while $\mathbf N$ is a commutative variety satisfying the identity~\eqref{xxy=0}.
\end{itemize}
\end{theorem}

Theorem~\ref{SEM umod deg>2} readily implies a necessary condition for a semigroup variety to be upper-modular given by~\cite[Theorem~1.1]{Vernikov-08-umod1} and the following

\begin{corollary}[\!\!{\mdseries\cite[Theorem~2]{Vernikov-Volkov-06}}]
\label{SEM umod nil}
A nil-variety of semigroups is upper-modular if and only if it satisfies the identities~\eqref{xxy=xyy,xy=yx}.
\end{corollary}

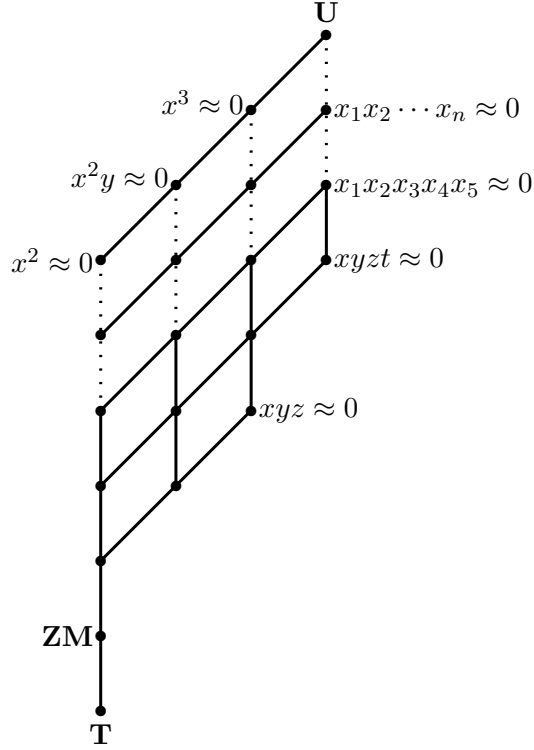
\begin{figure}[tbh]
\begin{center}
\unitlength=1mm
\linethickness{0.4pt}
\begin{picture}(70,100)
\put(12,5){\circle*{1.33}}
\put(12,15){\circle*{1.33}}
\put(12,25){\circle*{1.33}}
\put(12,35){\circle*{1.33}}
\put(12,45){\circle*{1.33}}
\put(12,55){\circle*{1.33}}
\put(12,65){\circle*{1.33}}
\put(22,35){\circle*{1.33}}
\put(22,45){\circle*{1.33}}
\put(22,55){\circle*{1.33}}
\put(22,65){\circle*{1.33}}
\put(22,75){\circle*{1.33}}
\put(32,45){\circle*{1.33}}
\put(32,55){\circle*{1.33}}
\put(32,65){\circle*{1.33}}
\put(32,75){\circle*{1.33}}
\put(32,85){\circle*{1.33}}
\put(42,65){\circle*{1.33}}
\put(42,75){\circle*{1.33}}
\put(42,85){\circle*{1.33}}
\put(42,95){\circle*{1.33}}
\gasset{AHnb=0,linewidth=0.4}
\drawline[dash={0.3 1.5}{0}](12,45)(12,65)
\drawline[dash={0.3 1.5}{0}](22,55)(22,75)
\drawline[dash={0.3 1.5}{0}](32,65)(32,85)
\drawline[dash={0.3 1.5}{0}](42,75)(42,95)
\drawline(12,5)(12,45)(42,75)(42,65)(32,55)(32,45)(12,25)
\drawline(12,35)(32,55)
\drawline(22,35)(22,55)
\drawline(12,55)(42,85)
\drawline(12,65)(42,95)
\drawline(32,55)(32,65)
\put(12,2){\makebox(0,0)[cc]{$\mathbf T$}}
\put(11,15){\makebox(0,0)[rc]{$\mathbf{ZM}$}}
\put(33,45){\makebox(0,0)[lc]{$xyz\approx0$}}
\put(43,65){\makebox(0,0)[lc]{$xyzt\approx0$}}
\put(43,75){\makebox(0,0)[lc]{$x_1x_2x_3x_4x_5\approx0$}}
\put(43,85){\makebox(0,0)[lc]{$x_1x_2\cdots x_n\approx0$}}
\put(11,65){\makebox(0,0)[rc]{$x^2\approx0$}}
\put(21,76){\makebox(0,0)[rc]{$x^2y\approx0$}}
\put(31,86){\makebox(0,0)[rc]{$x^3\approx0$}}
\put(42,98){\makebox(0,0)[cc]{$\mathbf U$}}
\end{picture}
\caption{The lattice $L(\mathbf U)$}
\label{SEM umod nil diagram}
\end{center}
\end{figure}

It is easy to list explicitly all upper-modular varieties of degree $>2$. Indeed, a simple arguments show that if a nil-variety satisfies the identity system~\eqref{xxy=xyy,xy=yx} then it satisfies also the identity $x^2yz\approx0$. Thus, it suffices to describe the subvariety lattice of the variety $\mathbf U=\var\{x^2y\approx xy^2,\,xy\approx yx,\,x^2yz\approx0\}$. Roughtine calculations show that this lattice has the form shown in Fig.\,\ref{SEM umod nil diagram}. The fact that some element of the lattice is marked in Fig.\,\ref{SEM umod nil diagram} by an identity $\varepsilon$ means that the corresponding variety is given within $\mathbf W$ by $\varepsilon$
. 

Theorem~\ref{SEM umod deg>2} reduces the examination of upper-modular varieties to varieties of degree $\le2$. To formulate a result concerning this case, we need some new definitions and notation. Recall that a semigroup variety is called \emph{completely regular} if it consists of \emph{completely regular} semigroups~--- unions of groups. A semigroup variety $\mathbf V$ is called a \emph{variety of semigroups with completely regular square} if, for any member $S$ of $\mathbf V$, the semigroup $S^2$ is completely regular. Put
\begin{align*}
&\mathbf{LZ}=\var\{xy\approx x\},\\
&\mathbf{RZ}=\var\{xy\approx y\},\\
&\mathbf P=\var\{xy\approx x^2y,\,x^2y^2\approx y^2x^2\},\\
&\overleftarrow{\mathbf P}=\var\{xy\approx xy^2,\,x^2y^2\approx y^2x^2\}.
\end{align*}
All we know about upper-modular varieties of degree $\le2$ is the following

\begin{theorem}[\!\!{\mdseries\cite[Theorem~2]{Vernikov-08-umod2}}]
\label{SEM umod deg<3-nec}
If $\mathbf V$ is an upper-modular semigroup variety of degree $\le2$ then one of the following holds:
\begin{itemize}
\item[\textup{(i)}]$\mathbf V$ is a variety of semigroups with completely regular square;
\item[\textup{(ii)}]$\mathbf{V=K\vee P}$ where $\mathbf K$ is a completely regular semigroup variety such that $\mathbf{RZ\nsubseteq K}$;
\item[\textup{(iii)}]$\mathbf{V=K}\vee\overleftarrow{\mathbf P}$ where $\mathbf K$ is a completely regular semigroup variety such that $\mathbf{LZ\nsubseteq K}$.
\end{itemize}
\end{theorem}

We do not know any example of a non-upper-modular variety that satisfies one of the claims~(i)--(iii) of Theorem~\ref{SEM umod deg<3-nec}. This inspires the following two questions.

\begin{question}
\label{SEM umod deg<3-crs?}
Is every variety of semigroups with completely regular square upper-modular?
\end{question}

\begin{question}
\label{SEM umod deg<3-lcr,rcr?}
Is every semigroup variety satisfying one of the claims~(ii) or~(iii) of Theorem~\ref{SEM umod deg<3-nec} upper-modular?
\end{question}

A natural weaker version of Question~\ref{SEM umod deg<3-crs?} is the following

\begin{question}
\label{SEM umod cr?}
Is every completely regular semigroup variety upper-modular?
\end{question}

Although Theorems~\ref{SEM umod deg>2} and~\ref{SEM umod deg<3-nec} do not provide a classification of all upper-modular varieties, they permit to deduct some important and surprising properties of such varieties. Theorems~\ref{SEM umod deg>2} and~\ref{SEM umod deg<3-nec} together with results of the articles~\cite{Volkov-92-mod3,Volkov-Ershova-91} imply the following

\begin{corollary}[\!\!{\mdseries\cite[Corollary~2]{Vernikov-08-umod2}}]
\label{SEM umod mod}
A proper upper-modular semigroup variety has a modular subvariety lattice.
\end{corollary}

Corollaries~\ref{SEM umod mod} and~\ref{SEM,COM umod,codistr ident&hered}(i) imply the following

\begin{corollary}[\!\!{\mdseries\cite[Corollary~3]{Vernikov-08-umod2}}]
\label{SEM umod hered}
If a proper semigroup variety is upper-modular then every its subvariety is also upper-modular.
\end{corollary}

Now we describe upper-modular varieties in one more class of varieties. A semigroup variety is called \emph{strongly permutative} if it satisfies an identity of the form $p_n[\pi]$ with $1\pi\ne1$ and $n\pi\ne n$.

\begin{theorem}
\label{SEM umod stpermut}
A strongly permutative semigroup variety $\mathbf V$ is upper-modular if and only if it satisfies one of the claims~\textup{(ii)} or~\textup{(iii)} of Theorem~\textup{\ref{SEM umod deg>2}}.
\end{theorem}

A partial case of this statement concerning commutative varieties is proved in~\cite[Theorem~1.2]{Vernikov-08-umod1}. Theorem~\ref{SEM umod stpermut} may be easily deduced from the proof of this partial case. A scheme of this deduction is provided in~\cite{Vernikov-08-lmod}.

As we have seen above (see Corollary~\ref{SEM umod mod}), the subvariety lattice of arbitrary proper upper-modular variety is modular. It turns out that such a lattice is even distributive in several wide classes of varieties. So, Theorem~\ref{SEM umod deg>2} together with results of the paper~\cite{Volkov-91} imply the following

\begin{corollary}[\!\!{\mdseries\cite[Corollary~1]{Vernikov-08-umod2}}]
\label{SEM umod deg>2-distr}
A proper upper-modular semigroup variety of degree $>2$ has a distributive subvariety lattice.
\end{corollary}

Theorem~\ref{SEM umod stpermut} together with results of~\cite{Volkov-91} imply the following

\begin{corollary}
\label{SEM umod stpermut-distr}
A strongly permutative upper-modular semigroup variety has a distributive subvariety lattice.
\end{corollary}

The special case of this claim dealing with commutative varieties was mentioned in~\cite[Corollary~4.4]{Vernikov-08-umod1}.

Theorem~\ref{SEM umod deg<3-nec} together with results of the articles~\cite{Rasin-82} and~\cite{Volkov-92-mod3} readily imply the following

\begin{corollary}[\!\!{\mdseries\cite[Corollary~4]{Vernikov-08-umod2}}]
\label{SEM umod not_crs-id}
Let $\mathbf V$ be a proper upper-modular semigroup variety that is not a variety of semigroups with completely regular square and let $\varepsilon$ be a non-trivial lattice identity. The lattice $L(\mathbf V)$ satisfies the identity $\varepsilon$ \textup(in particular, is distributive\textup) if and only if the subvariety lattice of any group subvariety of $\mathbf V$ has the same property.
\end{corollary}

Further, a semigroup variety $\mathbf V$ is called \emph{aperiodic} if all groups in $\mathbf V$ are trivial. Corollary~\ref{SEM umod not_crs-id} together with the result of the paper~\cite{Gerhard-77} readily imply the following

\begin{corollary}[\!\!{\mdseries\cite[Corollary~5]{Vernikov-08-umod2}}]
\label{SEM umod aper-distr}
An aperiodic upper-modular semigroup variety has a distributive subvariety lattice.
\end{corollary}

Corollaries~\ref{SEM umod deg>2-distr}--\ref{SEM umod aper-distr} inspire the following open

\begin{question}
\label{SEM umod distr?}
Is the subvariety lattice of every proper upper-modular semigroup variety distributive?
\end{question}

All proper upper-modular varieties that appeared above are varieties mentioned in Theorem~\ref{SEM umod deg>2}. These varieties are commutative. Based on this observation, one can conjecture that any proper upper-modular variety is commutative. But this is not the case. Evident counter-examples are the varieties $\mathbf{LZ}$ and $\mathbf{RZ}$. The claim that these two varieties are upper-modular immediately follows from the fact that they are atoms of the lattice $\mathbb{SEM}$. Two more examples of proper non-commutative upper-modular varieties are the varieties $\mathbf P$ and $\overleftarrow{\mathbf P}$. Indeed, it is well known that if a variety $\mathbf V$ is properly contained in one of these two varieties then $\mathbf{V\subseteq SL\vee ZM}$, whence $\mathbf V$ is lower-modular by Theorem~\ref{SEM lmod}. This readily implies that $\mathbf P$ and $\overleftarrow{\mathbf P}$ are upper-modular.

\subsection{Varieties that are both modular and upper-modular}
\label{subsection SEM:cumod}

It is interesting to examine varieties that satisfy different combinations of the properties we consider. Corollary~\ref{SEM lmod is cmod} implies that a variety is both modular and lower-modular if and only if it is lower-modular. So, Theorem~\ref{SEM lmod} gives, in fact, a complete description of varieties that are both modular and lower-modular (this result was obtained for the first time in~\cite[Theorem~3.1]{Volkov-05}). A description of varieties that both are lower-modular and upper-modular as well as varieties that are both distributive and codistributive is given in Theorem~\ref{SEM costand,neutr}. The following assertion classifies varieties that are both modular and upper-modular.

\begin{proposition}[\!\!{\mdseries\cite[Theorem~1]{Vernikov-Volkov-06}}]
\label{SEM cumod}
A semigroup variety $\mathbf V$ is both modular and upper-modular if and only if either $\mathbf{V=SEM}$ or $\mathbf{V=M\vee N}$ where $\mathbf M$ is one of the varieties $\mathbf T$ or $\mathbf{SL}$, while $\mathbf N$ is a commutative variety satisfying the identity~\eqref{xxy=0}.
\end{proposition}

\subsection{Codistributive varieties}
\label{subsection SEM:codistr}

The problem of description of codistributive semigroup varieties is open so far. Here we provide some partial results concerning this problem. The following theorem gives a strong necessary condition for a semigroup variety to be codistributive.

\begin{theorem}[\!\!{\mdseries\cite[Theorem~1.1]{Vernikov-11}}]
\label{SEM codistr nec}
If a semigroup variety $\mathbf V$ is codistributive then either $\mathbf{V=SEM}$ or $\mathbf V$ is a variety of semigroups with completely regular square.
\end{theorem}

Note that Theorems~\ref{SEM umod deg>2} and~\ref{SEM umod deg<3-nec} are crucial in the proof of Theorem~\ref{SEM codistr nec}.

It is easy to see that a variety of semigroups with completely regular square is a variety of degree $\le2$ (this readily follows from~\cite[Lemma~1]{Sapir-Sukhanov-81} or~\cite[Proposition~2.11]{Vernikov-08-umod1}). Therefore, Theorem~\ref{SEM codistr nec} implies that a proper codistributive variety has degree $\le2$. The following assertion shows that, for strongly permutative varieties, the converse statement holds as well.

\begin{theorem}[\!\!{\mdseries\cite[Theorem~1.2]{Vernikov-11}}]
\label{SEM codistr stpermut}
For a strongly permutative semigroup variety $\mathbf V$, the following are equivalent:
\begin{itemize}
\item[\textup{a)}]$\mathbf V$ is codistributive;
\item[\textup{b)}]$\mathbf V$ has a degree $\le2$;
\item[\textup{c)}]$\mathbf{V=A}_n\vee\mathbf X$ where $n\ge1$ and $\mathbf X$ is one of the varieties $\mathbf T$, $\mathbf{SL}$, $\mathbf{ZM}$ or $\mathbf{SL\vee ZM}$.
\end{itemize}
\end{theorem}

It is easy to see that there exist non-codistributive varieties of semigroups with completely regular square and moreover, non-codistributive periodic group varieties. Indeed, the lattice of periodic group varieties is modular but not distributive. Therefore, it contains the 5-element modular non-distributive sublattice. It is evident that all three pairwise non-comparable elements of this sublattice are non-codistributive periodic group varieties. We see that the problem of description of codistributive varieties is closely related to the problem of description of periodic group varieties with distributive subvariety lattice. The latter problem seems to be extremely difficult (see~\cite[Subsection~11.2]{Shevrin-Vernikov-Volkov-09} for more detailed comments), whence the former problem is extremely difficult too.

However, we do not know any examples of non-codistributive varieties of semigroups with completely regular square except ones mentioned in the previous paragraph. This inspires us to eliminate an examination of codistributive varieties with non-trivial groups. In other words, it seems natural to consider aperiodic codistributive varieties only. It is easy to see that if $\mathbf V$ is an aperiodic variety of semigroups with completely regular square then, for every $S\in\mathbf V$, the semigroups $S^2$ is a band. A variety with the last property is called a \emph{variety of semigroups with idempotent square}. In view of Theorem~\ref{SEM codistr nec}, every aperiodic codistributive variety is a variety of semigroups with idempotent square. Thus, the following question seems to be natural.

\begin{question}
\label{SEM codistr is?}
Is every variety of semigroups with idempotent square codistributive?
\end{question}

A natural weaker version of this question is the following

\begin{question}
\label{SEM codistr bands?}
Is every variety of bands codistributive?
\end{question}

Clearly, every variety of semigroups with idempotent square satisfies the identity $xy\approx(xy)^2$. Put
\begin{align*}
&\mathbf{IS}=\var\bigl\{xy\approx(xy)^2\bigr\},\\
&\mathbf{BAND}=\var\{x\approx x^2\}.
\end{align*}
It is verified in~\cite{Gerhard-77} that the lattice $L(\mathbf{IS})$ is distributive. Then Corollary~\ref{SEM,COM umod,codistr ident&hered}(ii) shows that Question~\ref{SEM codistr is?} is equivalent to the following: is the variety $\mathbf{IS}$ codistributive? Similarly, Question~\ref{SEM codistr bands?} is equivalent to asking, whether the variety $\mathbf{BAND}$ is codistributive, i.e., whether the equality
$$
\mathbf{BAND\wedge(X\vee Y)=(BAND\wedge X)\vee(BAND\wedge Y)}
$$
holds for arbitrary varieties $\mathbf X$ and $\mathbf Y$ or not. It is verified in~\cite[Corollary~5.9]{Pastijn-Trotter-98} that this is the case whenever the varieties $\mathbf X$ and $\mathbf Y$ are locally finite.

A strongly permutative codistributive variety has a distributive subvariety lattice (this follows from Corollary~\ref{SEM umod stpermut-distr} and may be easily deduced from Theorem~\ref{SEM codistr stpermut}). Aperiodic codistributive varieties also have a distributive subvariety lattice (here it suffices to refer to either Corollary~\ref{SEM umod aper-distr} or Theorem~\ref{SEM codistr nec} and the mentioned result of~\cite{Gerhard-77}). We do not know any example of proper codistributive variety with non-distributive subvariety lattice. This inspires the following

\begin{question}
\label{SEM codistr distr?}
Is the subvariety lattice of every proper codistributive semigroup variety distributive?
\end{question}

This question is closely related to the following

\begin{question}
\label{SEM codistr hered?}
Is every subvariety of an arbitrary proper codistributive semigroup variety codistributive?
\end{question}

Corollary~\ref{SEM,COM umod,codistr ident&hered}(ii) shows that the affirmative answer to Question~\ref{SEM codistr distr?} would imply the affirmative answer to Question~\ref{SEM codistr hered?}.

As we have mentioned in Introduction, results of the article~\cite{Aizenstat-74} imply that all atoms of the lattice $\mathbb{SEM}$ are codistributive. This fact can be generalized by the following way.

\begin{remark}[\!\!{\mdseries\cite[Remark~4.1]{Vernikov-11}}]
\label{SEM codistr fin-join-of-atoms}
If $\mathbf V_1$, $\mathbf V_2$, \dots, $\mathbf V_k$ are atoms of the lattice $\mathbb{SEM}$ then the variety $\mathop\bigvee\limits_{i=1}^k\mathbf V_i$ is codistributive.
\end{remark}

In connection with Questions~\ref{SEM codistr distr?} and~\ref{SEM codistr hered?}, we note that if $\mathbf V_1$, $\mathbf V_2$, \dots, $\mathbf V_k$ are atoms of the lattice $\mathbb{SEM}$ and $\mathbf V=\mathop\bigvee\limits_{i=1}^k\mathbf V_i$ then:
\begin{itemize}
\item[(i)]the lattice $L(\mathbf V)$ is distributive (in fact, $L(\mathbf V)$ is a direct product of $k$ copies of the 2-element chain),
\item[(ii)]if $\mathbf{X\subseteq V}$ then $\mathbf X$ is the join of those of the atoms $\mathbf V_1$, $\mathbf V_2$, \dots, $\mathbf V_k$ that are contained in $\mathbf X$, and therefore, $\mathbf X$ is codistributive by Remark~\ref{SEM codistr fin-join-of-atoms}.
\end{itemize}
The claim~(i) is a part of~\cite[Proposition~1]{Vernikov-Volkov-82}, while the statement~(ii) follows from~(i).

\section{The lattice $\mathbb{COM}$}
\label{section COM}

For convenience, we call a commutative semigroup variety $\mathbb{COM}$-\emph{modular} if it is a modular element of the lattice $\mathbb{COM}$ and adopt similar agreement for all other types of special elements. The main results of this section provide:
\begin{itemize}
\item a complete classification of $\mathbb{COM}$-lower-modular, $\mathbb{COM}$-upper-modular, $\mathbb{COM}$-distributive, $\mathbb{COM}$-codistributive, $\mathbb{COM}$-standard, $\mathbb{COM}$-costandard or $\mathbb{COM}$-neutral varieties (Theorems~\ref{COM lmod},~\ref{COM distr,stand},~\ref{COM neutr},~\ref{COM umod,codistr} and~\ref{COM costand}),{\sloppy

}
\item necessary conditions for a commutative semigroup variety to be $\mathbb{COM}$-modular (Theorems~\ref{COM cmod nec} and~\ref{COM cmod nil-nec}),
\item a sufficient condition for a commutative semigroup variety to be $\mathbb{COM}$-modular (Theorem~\ref{COM cmod suf}).
\end{itemize}

Note that there are no any essential information about $\mathbb{COM}$-cancellable varieties so far.

\subsection{$\mathbb{COM}$-lower-modular varieties}
\label{subsection COM:lmod}

We denote by $\mathbf{COM}$ the variety of all commutative semigroups. A commutative semigroup variety is called $\mathbb{COM}$-0-\emph{reduced} if it may be given by the commutative law and some non-empty set of 0-reduced identities only. Some partial information about $\mathbb{COM}$-lower-modular varieties was obtained in~\cite{Shaprynskii-11}. It is covered by the following `commutative analogue' of Theorem~\ref{SEM lmod}.

\begin{theorem}[\!\!{\mdseries\cite[Theorem~1.6]{Shaprynskii-12-mod-lmod}}]
\label{COM lmod}
A commutative semigroup variety $\mathbf V$ is $\mathbb{COM}$-lower-modular if and only if either $\mathbf{V=COM}$ or $\mathbf{V=M\vee N}$ where $\mathbf M$ is one of the varieties $\mathbf T$ or $\mathbf{SL}$, while $\mathbf N$ is a $\mathbb{COM}$-$0$-reduced variety.
\end{theorem}

Note that the `if' part of Theorem~\ref{COM lmod} immediately follows from Corollaries~\ref{L(X) cmod,lmod suf} (with $\mathbf{X=COM}$) and~\ref{join with SL I-elem}. The proof of the `only if' part given in~\cite{Shaprynskii-12-mod-lmod} is based on Theorem~\ref{L(oc) cmod,lmod per-nec} below.

As in the case of the lattice $\mathbb{SEM}$ (see Subsection~\ref{subsection SEM:lmod}), Theorem~\ref{COM lmod} implies that a description of $\mathbb{COM}$-distributive, $\mathbb{COM}$-standard and $\mathbb{COM}$-neutral varieties should look as follows:

\emph{A commutative semigroup variety $\mathbf V$ is $\mathbb{COM}$-distributive \textup[$\mathbb{COM}$-standard, $\mathbb{COM}$-neutral\textup] if and only if either $\mathbf{V=COM}$ or $\mathbf{V=M\vee N}$ where $\mathbf M$ is one of the varieties $\mathbf T$ or $\mathbf{SL}$, while $\mathbf N$ is a $\mathbb{COM}$-$0$-reduced variety such that \dots} (with some additional restriction to $\mathbf N$ depending on the type of element we consider).

Exact formulations of corresponding results are given in the following two subsections.

\subsection{$\mathbb{COM}$-distributive and $\mathbb{COM}$-standard varieties}
\label{subsection COM:distr,stand}

The following statement is the `commutative analogue' of Theorem~\ref{SEM distr,stand}.

\begin{theorem}[\!\!{\mdseries\cite[Theorem~1.1]{Shaprynskii-11}}]
\label{COM distr,stand}
For a commutative semigroup variety $\mathbf V$, the following are equivalent:
\begin{itemize}
\item[\textup{a)}]$\mathbf V$ is $\mathbb{COM}$-distributive;
\item[\textup{b)}]$\mathbf V$ is $\mathbb{COM}$-standard;
\item[\textup{c)}]either $\mathbf{V=COM}$ or $\mathbf{V=M\vee N}$ where $\mathbf M$ is one of the varieties $\mathbf T$ or $\mathbf{SL}$, while $\mathbf N$ is a $\mathbb{COM}$-$0$-reduced variety that satisfies the identities $x^3yz\approx x^2y^2z\approx0$ and either satisfies both the identities $x^3y\approx0$ and $x^2y^2\approx0$ or does not satisfy any of them.
\end{itemize}
\end{theorem}

Note that the item c) of Theorem \ref{COM distr,stand} is reproduced in~\cite{Vernikov-15} with inaccuracy\footnote{Namely, the identity $x^3\approx0$ rather than $x^3y\approx0$ is written in the item c) of~\cite[Theorem~4.2]{Vernikov-15}.}.

It is verified in~\cite[Corollary~1.1]{Shaprynskii-11} that a $\mathbb{COM}$-distributive variety is $\mathbb{COM}$-modular. This statement is generalized by Corollary~\ref{COM lmod is cmod} below.

\subsection{$\mathbb{COM}$-neutral varieties}
\label{subsection COM:neutr}

A complete description of $\mathbb{COM}$-neutral varieties is given by the following partial analogue of Theorem~\ref{SEM costand,neutr}.

\begin{theorem}
\label{COM neutr}
For a commutative semigroup variety $\mathbf V$, the following are equivalent:
\begin{itemize}
\item[\textup{a)}]$\mathbf V$ is both $\mathbb{COM}$-upper-modular and $\mathbb{COM}$-lower-modular;
\item[\textup{b)}]$\mathbf V$ is both $\mathbb{COM}$-distributive and $\mathbb{COM}$-codistributive;
\item[\textup{c)}]$\mathbf V$ is $\mathbb{COM}$-neutral;
\item[\textup{d)}]either $\mathbf{V=COM}$ or $\mathbf{V=M\vee N}$ where $\mathbf M$ is one of the varieties $\mathbf T$ or $\mathbf{SL}$ and the variety $\mathbf N$ satisfies the identity~\eqref{xxy=0}.
\end{itemize}
\end{theorem}

The equivalence of the claims~b)--d) of this theorem is verified in~\cite[Theorem~1.2]{Shaprynskii-11}, while the equivalence~a)\,$\Longleftrightarrow$\,c) is proved in~\cite[Corollary~4.2]{Shaprynskii-12-mod-lmod}.

Theorem~\ref{COM neutr} and Corollary~\ref{SEM cmod,canc commut} immediately implies the following

\begin{corollary}
\label{COM neutr = SEM cmod,canc}
For a commutative semigroup variety $\mathbf V$ with $\mathbf{V\ne COM}$, the following are equivalent:
\begin{itemize}
\item[\textup{a)}]$\mathbf V$ is $\mathbb{COM}$-neutral;
\item[\textup{b)}]$\mathbf V$ is modular;
\item[\textup{c)}]$\mathbf V$ is cancellable.
\end{itemize}
\end{corollary}

\subsection{$\mathbb{COM}$-modular varieties}
\label{subsection COM:cmod} 

The problem of description of $\mathbb{COM}$-modular commutative semigroup varieties is open so far. Here we provide some partial results concerning this problem. Note that these results are `commutative analogues' of Theorems~\ref{SEM cmod nec},~\ref{SEM cmod nil-nec} and~\ref{SEM cmod suf}.

First of all, the following necessary condition for a commutative semigroup variety to be $\mathbb{COM}$-modular is true.

\begin{theorem}[\!\!{\mdseries\cite[Theorem~1.4]{Shaprynskii-12-mod-lmod}}]
\label{COM cmod nec}
If $\mathbf V$ is a $\mathbb{COM}$-modular commutative semigroup variety then either $\mathbf{V=COM}$ or $\mathbf{V=M\vee N}$ where $\mathbf M$ is one of the varieties $\mathbf T$ or $\mathbf{SL}$, while $\mathbf N$ is a nil-variety.
\end{theorem}

In fact, this theorem readily follows from Theorem~\ref{L(oc) cmod,lmod per-nec} below. Theorem~\ref{COM cmod nec} and Corollary~\ref{join with SL I-elem} completely reduce an examination of $\mathbb{COM}$-modular varieties to the nil-case. The following theorem is yet another analogue of the result of Je\v{z}ek~\cite{Jezek-81} (see Theorem~\ref{SEM cmod nil-nec} and the paragraph before this theorem).

\begin{theorem}[\!\!{\mdseries\cite[Theorem~1.5]{Shaprynskii-12-mod-lmod}}]
\label{COM cmod nil-nec}
A $\mathbb{COM}$-modular commutative nil-variety of semigroups may be given within the variety $\mathbf{COM}$ by $0$-reduced and substitutive identities only.
\end{theorem}

Corollary~\ref{L(X) cmod,lmod suf} with $\mathbf{X=COM}$ immediately implies the following

\begin{theorem}[\!\!{\mdseries\cite[Proposition~2.1]{Shaprynskii-11}}]
\label{COM cmod suf}
Every $\mathbb{COM}$-$0$-reduced commutative semigroup variety is $\mathbb{COM}$-modular.
\end{theorem}

Theorems~\ref{COM cmod nil-nec} and~\ref{COM cmod suf} provide, respectively, a necessary and a sufficient condition for a commutative nil-variety to be $\mathbb{COM}$-modular. The gap between these conditions does not seem to be very large. But the necessary condition is not a sufficient one, while the sufficient condition is not a necessary one. Indeed:
\begin{itemize}
\item the variety $\var\{xyzt\approx x^3\approx0,\,x^2y\approx y^2x,\,xy\approx yx\}$ is $\mathbb{COM}$-modular although it is not $\mathbb{COM}$-0-reduced~\cite[Proposition~4.1]{Shaprynskii-15},
\item the variety\\
\rule{3mm}{0pt}$\var\{x^4y^3z^2t\approx y^4x^3z^2t,\,x_1x_2\cdots x_{11}\approx\mathbf w_i\approx 0,\,xy\approx yx\mid\mathbf w_i\in W_{10,3}\}$\\
where $W_{10,3}$ is the set of all words of length 10 depending on $\le3$ letters is not $\mathbb{COM}$-modular although it is given within $\mathbb{COM}$ by 0-reduced and substitutive identities only~\cite[Proposition~4.2]{Shaprynskii-15}.
\end{itemize}

Theorems~\ref{COM lmod} and~\ref{COM cmod suf} and Corollary~\ref{join with SL I-elem} imply the following `commutative analogue' of Corollary~\ref{SEM lmod is cmod}.

\begin{corollary}[\!\!{\mdseries\cite[Corollary~4.1]{Shaprynskii-12-mod-lmod}}]
\label{COM lmod is cmod}
Every $\mathbb{COM}$-lower-modular commutative semigroup variety is $\mathbb{COM}$-modular.
\end{corollary}

\subsection{$\mathbb{COM}$-cancellable varieties}
\label{subsection COM:canc}

There are no any significant information about $\mathbb{COM}$-cancellable commutative semigroup varieties so far. We mention only a few observations here. 

First, it is clear that the conclusions of Theorems~\ref{COM cmod nec} and~\ref{COM cmod nil-nec} remain true whenever $\mathbf V$ is a $\mathbb{COM}$-cancellable commutative semigroup variety.

Second, the `cancellable analogue' of Theorem~\ref{COM cmod suf} is false. This immediately follows from the complete description of the subvariety lattice of the variety
$$
\mathbf N_5^c=\var\{x_1x_2x_3x_4x_5\approx0,\,xy\approx yx\}
$$
obtained in~\cite[Fig.~3]{Melnik-72}. We reproduce the diagram of this lattice in Fig.~\ref{L(N_5^c)} below. This figure shows that the $\mathbb{COM}$-0-reduced varieties $\var\{x_1x_2x_3x_4x_5\approx x^3\approx0,\,xy\approx yx\}$ and $\var\{x_1x_2x_3x_4x_5\approx x^3y\approx0,\,xy\approx yx\}$ are not $\mathbb{COM}$-cancellable. 

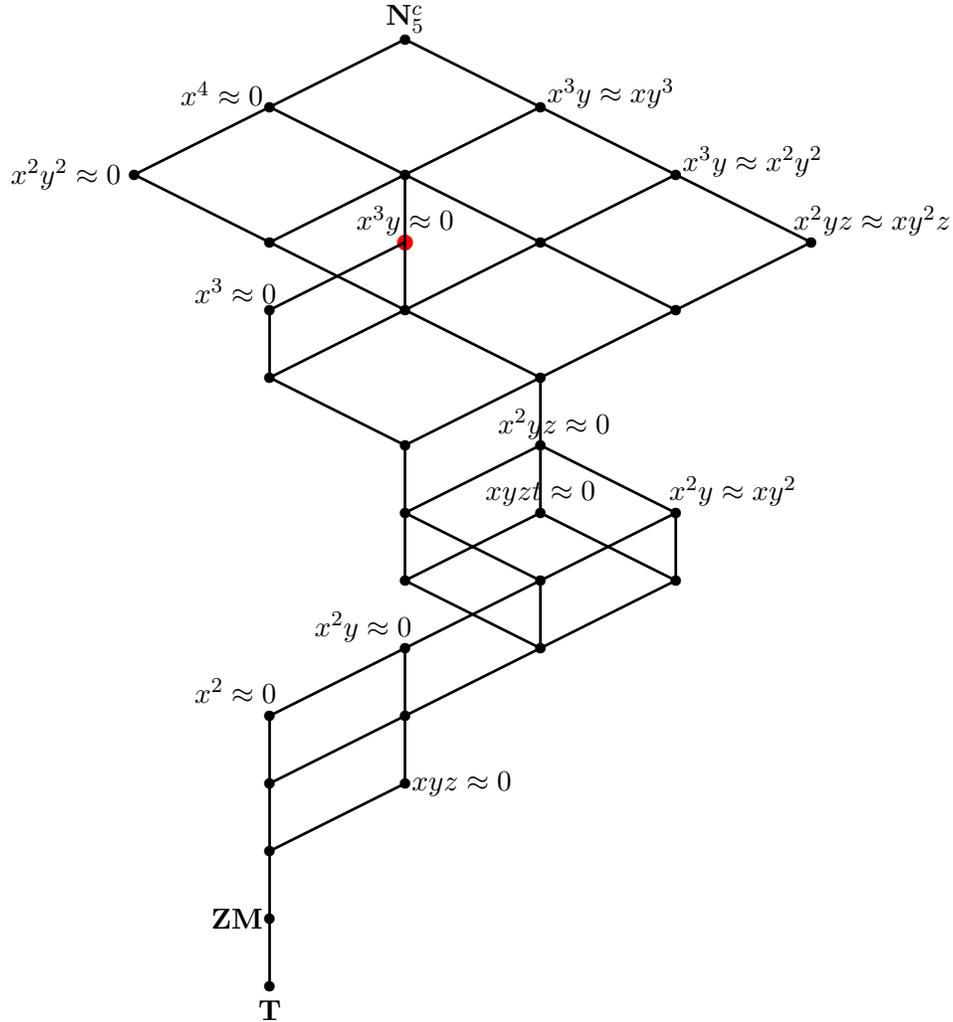
\begin{figure}[htb]
\unitlength=1mm
\linethickness{0.4pt}
\begin{center}
\gasset{AHnb=0,linewidth=0.4}
\begin{picture}(132,152)
\put(32,6){\circle*{1.33}}
\put(32,16){\circle*{1.33}}
\put(32,26){\circle*{1.33}}
\put(32,36){\circle*{1.33}}
\put(52,36){\circle*{1.33}}
\put(32,46){\circle*{1.33}}
\put(52,46){\circle*{1.33}}
\put(52,56){\circle*{1.33}}
\put(72,56){\circle*{1.33}}
\put(52,66){\circle*{1.33}}
\put(72,66){\circle*{1.33}}
\put(92,66){\circle*{1.33}}
\put(52,76){\circle*{1.33}}
\put(72,76){\circle*{1.33}}
\put(92,76){\circle*{1.33}}
\put(52,86){\circle*{1.33}}
\put(72,86){\circle*{1.33}}
\put(32,96){\circle*{1.33}}
\put(72,96){\circle*{1.33}}
\put(32,106){\circle*{1.33}}
\put(52,106){\circle*{1.33}}
\put(92,106){\circle*{1.33}}
\put(32,116){\circle*{1.33}}
\put(52,116){\circle*{1.33}}
\put(72,116){\circle*{1.33}}
\put(112,116){\circle*{1.33}}
\put(12,126){\circle*{1.33}}
\put(52,126){\circle*{1.33}}
\put(92,126){\circle*{1.33}}
\put(32,136){\circle*{1.33}}
\put(72,136){\circle*{1.33}}
\put(52,146){\circle*{1.33}}
\drawline(32,6)(32,46)(92,76)(92,66)(32,36)
\drawline(32,26)(52,36)(52,56)
\drawline(52,66)(52,86)(112,116)(52,146)(12,126)(72,96)(72,76)(52,66)(72,56)(72,66)(52,76)(72,86)(92,76)
\drawline(72,76)(92,66)
\drawline(52,86)(32,96)(92,126)
\drawline(32,96)(32,106)(52,116)(52,106)
\drawline(52,116)(52,126)(32,116)
\drawline(32,136)(92,106)
\drawline(52,126)(72,136)
\put(32,4){\makebox(0,0)[ct]{$\mathbf T$}}
\put(31,16){\makebox(0,0)[rc]{$\mathbf{ZM}$}}
\put(52,147){\makebox(0,0)[cb]{$\mathbf N_5^c$}}
\put(8,123){\makebox(0,0)[cc]{$x^2y^2\approx0$}}
\put(33,107){\makebox(0,0)[rb]{$x^3\approx0$}}
\put(31,138){\makebox(0,0)[rc]{$x^4\approx0$}}
\put(33,48){\makebox(0,0)[rb]{$x^2\approx0$}}
\put(52,117){\makebox(0,0)[cb]{$x^3y\approx0$}}
\put(53,36){\makebox(0,0)[lc]{$xyz\approx0$}}
\put(53,57){\makebox(0,0)[rb]{$x^2y\approx0$}}
\put(72,77){\makebox(0,0)[cb]{$xyzt\approx0$}}
\put(73,138){\makebox(0,0)[lc]{$x^3y\approx xy^3$}}
\put(74,87){\makebox(0,0)[cb]{$x^2yz\approx0$}}
\put(91,77){\makebox(0,0)[lb]{$x^2y\approx xy^2$}}
\put(93,128){\makebox(0,0)[lc]{$x^3y\approx x^2y^2$}}
\put(121,113){\makebox(0,0)[cc]{$x^2yz\approx xy^2z$}}
\end{picture}
\end{center}
\caption{The lattice $L(\mathbf N_5^c)$}
\label{L(N_5^c)}
\end{figure}

\subsection{$\mathbb{COM}$-upper-modular and $\mathbb{COM}$-codistributive varieties}
\label{subsection COM:umod,codistr}

A complete description of $\mathbb{COM}$-upper-modular and $\mathbb{COM}$-codistributive varieties is given by the following{\sloppy

}

\begin{theorem}[\!\!{\mdseries\cite[Theorem~1.1]{Vernikov-17}}]
\label{COM umod,codistr}
For a commutative semigroup variety $\mathbf V$, the following are equivalent:
\begin{itemize}
\item[\textup{a)}]$\mathbf V$ is $\mathbb{COM}$-upper-modular;
\item[\textup{b)}]$\mathbf V$ is $\mathbb{COM}$-codistributive;
\item[\textup{c)}]either $\mathbf{V=COM}$ or $\mathbf V$ satisfies one of the claims~\textup{(ii)} or~\textup{(iii)} of Theorem~\textup{\ref{SEM umod deg>2}}.
\end{itemize}
\end{theorem}

Note that, in contrast with the equivalence of the claims~a) and~b) of this theorem, the properties of being upper-modular and codistributive varieties are not equivalent. This follows from comparison of Theorems~\ref{SEM umod stpermut} and~\ref{SEM codistr stpermut}.

Theorems~\ref{SEM umod stpermut} and~\ref{COM umod,codistr} implies immediately

\begin{corollary}[\!\!{\mdseries\cite[Corollary~4.1]{Vernikov-17}}]
\label{umod COM and SEM}
If $\mathbf V$ is a commutative semigroup variety and $\mathbf{V\ne COM}$ then the following are equivalent:
\begin{itemize}
\item[\textup{a)}]$\mathbf V$ is $\mathbb{COM}$-upper-modular;
\item[\textup{b)}]$\mathbf V$ is $\mathbb{COM}$-codistributive;
\item[\textup{c)}]$\mathbf V$ is upper-modular.
\end{itemize}
\end{corollary}

Theorem~\ref{COM neutr} implies the following `commutative analogue' of Corollary~\ref{SEM umod hered}.

\begin{corollary}[\!\!{\mdseries\cite[Corollary~4.3]{Vernikov-17}}]
\label{COM umod hered}
If $\mathbf V$ is a $\mathbb{COM}$-upper-modular commutative semigroup variety and $\mathbf{V\ne COM}$ then every subvariety of $\mathbf V$ is $\mathbb{COM}$-upper-modular.
\end{corollary}

Theorem~\ref{COM umod,codistr} together with results of~\cite{Volkov-91} imply the following

\begin{corollary}[\!\!{\mdseries\cite[Corollary~4.4]{Vernikov-17}}]
\label{COM umod distr}
If $\mathbf V$ is a $\mathbb{COM}$-upper-modular commutative semigroup variety and $\mathbf{V\ne COM}$ then the lattice $L(\mathbf V)$ is distributive.
\end{corollary}

\subsection{$\mathbb{COM}$-costandard varieties}
\label{subsection COM:costand}

A complete description of $\mathbb{COM}$-costandard varieties is given by the following

\begin{theorem}[\!\!{\mdseries\cite[Theorem~1.2]{Vernikov-17}}]
\label{COM costand}
For a commutative semigroup variety $\mathbf V$, the following are equivalent:
\begin{itemize}
\item[\textup{a)}]$\mathbf V$ is both $\mathbb{COM}$-modular and $\mathbb{COM}$-upper-modular;
\item[\textup{b)}]$\mathbf V$ is $\mathbb{COM}$-costandard;
\item[\textup{c)}]either $\mathbf{V=COM}$ or $\mathbf V$ satisfies the claim~\textup{(ii)} of Theorem~\textup{\ref{SEM umod deg>2}}.
\end{itemize}
\end{theorem}

A comparison of Theorems~\ref{COM neutr} and~\ref{COM costand} shows that, in contrast with Theorem~\ref{SEM costand,neutr}, the properties of being $\mathbb{COM}$-neutral and $\mathbb{COM}$-costandard varieties are not equivalent.

As we have already mentioned in Subsection \ref{subsection SEM:umod}, in every nil-variety, the identities \eqref{xxy=xyy,xy=yx} imply the identity $x^2yz\approx0$. Then comparison of Theorems \ref{COM distr,stand} and \ref{COM costand} implies the following analog of Corollary \ref{SEM costand is stand} which was not mentioned explicitly anywhere.

\begin{corollary}
\label{COM costand is stand}
Every $\mathbb{COM}$-costandard commutative semigroup variety is $\mathbb{COM}$-standard.
\end{corollary}

Note that a description of varieties that are both $\mathbb{COM}$-modular and $\mathbb{COM}$-lower-modular immediately follows from Corollary~\ref{COM lmod is cmod} and Theorem~\ref{COM lmod}, while a description of varieties that are both $\mathbb{COM}$-upper-modular and $\mathbb{COM}$-lower-modular is given by Theorem~\ref{COM neutr}. We note also that the analogue of Theorem~\ref{COM costand} in the lattice $\mathbb{SEM}$ is not the case. Indeed, a comparison of Theorem~\ref{SEM costand,neutr} and Proposition~\ref{SEM cumod} show that there exist semigroup varieties that are simultaneously modular and upper-modular but not costandard.

Theorems~\ref{COM umod,codistr} and~\ref{COM costand} imply immediately

\begin{corollary}[\!\!{\mdseries\cite[Corollary~4.2]{Vernikov-17}}]
\label{COM umod,codistr,costand nil}
For a commutative nil-variety of semigroups $\mathbf V$, the following are equivalent:
\begin{itemize}
\item[\textup{a)}]$\mathbf V$ is $\mathbb{COM}$-upper-modular;
\item[\textup{b)}]$\mathbf V$ is $\mathbb{COM}$-codistributive;
\item[\textup{c)}]$\mathbf V$ is $\mathbb{COM}$-costandard;
\item[\textup{d)}]$\mathbf V$ satisfies the identities~\eqref{xxy=xyy,xy=yx}.
\end{itemize}
\end{corollary}

Note that the equivalence of the claims~a)--c) of Corollary~\ref{COM umod,codistr,costand nil} for $\mathbb{COM}$-0-reduced varieties immediately follows from~\cite[Theorem~1.2 and Proposition~2.3]{Shaprynskii-11}. Moreover, if a variety $\mathbf V$ is $\mathbb{COM}$-0-reduced then the claims~a)--c) are equivalent to the claim that $\mathbf V$ is $\mathbb{COM}$-neutral. This also follows from the mentioned results of~\cite{Shaprynskii-11}.

\section{Lattices located between $\mathbb{SEM}$ and $\mathbb{COM}$}
\label{section between SEM and COM}

In this section, we examine modular and lower-modular elements only. It turns out that properties of such elements in the lattices $\mathbb{SEM}$ and $\mathbb{COM}$ discussed in Subsections~\ref{subsection SEM:lmod},~\ref{subsection SEM:cmod},~\ref{subsection COM:lmod} and~\ref{subsection COM:cmod} may be partially extended to some sublattices of $\mathbb{SEM}$ that contain $\mathbb{COM}$. More precisely, we have in mind subvariety lattices of overcommutative semigroup varieties and the lattice $\mathbb{PERM}$.

\subsection{Subvariety lattices of overcommutative varieties}
\label{subsection between-SEM-and-COM:L(oc)}

As we have seen in Subsections~\ref{subsection COM:lmod} and~\ref{subsection COM:cmod}, there are numerous parallels between results about modular and lower-modular elements in the lattices~$\mathbb{SEM}$ and~$\mathbb{COM}$. The following result partially explains these parallels and permits us to give unified proofs of several results about [lower-]modular elements in $\mathbb{SEM}$ and $\mathbb{COM}$.

\begin{theorem}[\!\!{\mdseries\cite[Proposition~3.3]{Shaprynskii-12-mod-lmod}}]
\label{L(oc) cmod,lmod per-nec}
Let $\mathbf X$ be an overcommutative semigroup variety and $\mathbf V$ a periodic subvariety of $\mathbf X$. If $\mathbf V$ is either a modular or a lower-modular element of the lattice $L(\mathbf X)$ then $\mathbf{V=M\vee N}$ where $\mathbf M$ is one of the varieties $\mathbf T$ or $\mathbf{SL}$, while $\mathbf N$ is a nil-variety.
\end{theorem}

Applying this theorem with $\mathbf{X=SEM}$ [respectively $\mathbf{X=COM}$], we obtain an important information about [$\mathbb{COM}$-]modular and [$\mathbb{COM}$-]lower-modular varieties. After that, only some simple additional arguments are needed to verify Theorems~\ref{SEM cmod suf} and~\ref{COM cmod suf}, as well as the `if' parts of Theorems~\ref{SEM lmod} and~\ref{COM lmod}. It is natural to ask if it is possible to eliminate these additional arguments altogether. To do this, we should verify an analogue of Theorem~\ref{L(oc) cmod,lmod per-nec} without the assumption that the variety $\mathbf V$ is periodic. Unfortunately, it turns out that this is impossible. Indeed, it is verified in~\cite{Trakhtman-74} that every proper semigroup variety is covered in $\mathbb{SEM}$ by some other variety (see also~\cite[Subsection~3.1]{Shevrin-Vernikov-Volkov-09}). It is evident that if an overcommutative variety $\mathbf V$ is covered by a variety $\mathbf X$ then $\mathbf X$ is overcommutative and $\mathbf V$ is a lower-modular element of the lattice $L(\mathbf X)$. Thus, the `lower-modular half' of Theorem~\ref{L(oc) cmod,lmod per-nec} would be false if we eliminate the assumption that $\mathbf V$ is periodic. The same is true for the `modular half' of this theorem. For example, the variety $\mathbf{COM}$ is a modular element in the lattice $L(\mathbf W)$ where $\mathbf W=\var\{xyz\approx yzx\approx zyx\}$~\cite[p.\,29]{Shaprynskii-15}. Note that $\mathbf{COM}$ is also a lower-modular element in $L(\mathbf W)$ because $\mathbf W$ covers $\mathbf{COM}$.

\subsection{The lattice $\mathbb{PERM}$}
\label{subsection between-SEM-and-COM:PERM}

By analogy with the commutative case, we call a permutative semigroup variety $\mathbb{PERM}$-[\emph{lower-}]\emph{modular} if it is a [lower-]modular element of the lattice $\mathbb{PERM}$.

\begin{theorem}[\!\!{\mdseries\cite[Proposition~2.2]{Shaprynskii-15}}]
\label{PERM cmod,lmod nec}
If a permutative semigroup variety $\mathbf V$ is either $\mathbb{PERM}$-modular or $\mathbb{PERM}$-lower-modular then $\mathbf{V=M\vee N}$ where $\mathbf M$ is one of the varieties $\mathbf T$ or $\mathbf{SL}$, while $\mathbf N$ is a nil-variety.
\end{theorem}

This result does not give any information about $\mathbb{PERM}$-modular or $\mathbb{PERM}$-lower-modular nil-varieties. Recall that:
\begin{itemize}
\item[(i)]by Theorems~\ref{SEM cmod nil-nec} and~\ref{COM cmod nil-nec}, every [$\mathbb{COM}$-]modular nil-variety may be given [within $\mathbf{COM}$] by substitutive and 0-reduced identities only;
\item[(ii)]by Theorems~\ref{SEM lmod} and~\ref{COM lmod}, every [$\mathbb{COM}$-]lower-modular nil-variety is [$\mathbb{COM}$-]0-reduced;
\item[(iii)]by Corollary~\ref{L(X) cmod,lmod suf}, every [$\mathbb{COM}$-]0-reduced variety is both [$\mathbb{COM}$-]modular and [$\mathbb{COM}$-]lower-modular.
\end{itemize}
Note that we cannot use Corollary~\ref{L(X) cmod,lmod suf} to obtain a `permutative analogue' of the claim~(iii) because the class of all permutative semigroups does not form a variety.

We do not know, whether a `permutative analogue' of the claim~(i) true. So, we formulate the following

\begin{question}
\label{PERM cmod nil-nec?}
Is it true that every $\mathbb{PERM}$-modular permutative nil-variety of semigroups may be given by substitutive and 0-reduced identities only?
\end{question}

As to `permutative analogues' of claims~(ii) and~(iii), they do not hold. For instance:
\begin{itemize}
\item the variety $\var\{xyzt\approx0,\,x^2y\approx xyx\}$ is $\mathbb{PERM}$-lower-modular although it may not be given by permutational and 0-reduced identities only;
\item the variety $\mathbf N_5^c $ is neither $\mathbb{PERM}$-modular nor $\mathbb{PERM}$-lower-modular although it is permutative and is given by permutational and 0-reduced identities only (the variety $\mathbf N_5^c $ was introduced in Subsection~\ref{subsection COM:canc}).
\end{itemize}
Both these claims are communicated to the author by Shaprynski\v{\i}.

\section{The lattice $\mathbb{OC}$}
\label{section OC}

For convenience, we call an overcommutative semigroup variety $\mathbb{OC}$-\emph{modular} if it is a modular element of the lattice $\mathbb{OC}$ and adopt similar agreement for all other types of special elements.

The problems of description of $\mathbb{OC}$-modular, $\mathbb{OC}$-lower-modular and $\mathbb{OC}$-upper-modular varieties are open so far. Moreover, any essential information about varieties of these three types is absent. On the other hand, $\mathbb{OC}$-distributive, $\mathbb{OC}$-codistributive, $\mathbb{OC}$-cancellable, $\mathbb{OC}$-standard, $\mathbb{OC}$-costandard and $\mathbb{OC}$-neutral varieties are completely determined (see Theorem~\ref{OC neutr,stand,costand,canc,distr,codistr} below). To formulate this result, we need some new definitions and notation.

Let $m$ and $n$ be positive integers with $2\le m\le n$. A sequence of positive integers $(\ell_1,\ell_2,\dots,\ell_m)$ is called a \emph{partition of $n$ into $m$ parts} if
$$
\sum_{i=1}^m\ell_i=n\quad\text{and}\quad\ell_1\ge\ell_2\ge\cdots\ge\ell_m.
$$
The set of all partitions of $n$ into $m$ parts is denoted by $\Lambda_{n,m}$. Let $\lambda=(\ell_1,\ell_2,\dots,\ell_m)\in\Lambda_{n,m}$. We define numbers $q(\lambda)$, $r(\lambda)$ and $s(\lambda)$ as follows:

$q(\lambda)$ is the number of $\ell_i$'s with $\ell_i=1$;

$r(\lambda)=n-q(\lambda)$ (in other words, $r(\lambda)$ is the sum of all $\ell_i$'s with $\ell_i>1$);

$s(\lambda)=\max\bigl\{r(\lambda)-q(\lambda)-\delta,0\bigr\}$ where
$$
\delta=
\begin{cases}
0&\text{if } n=3,\,m=2\text{ and }\lambda=(2,1),\\
1&\text{otherwise}.
\end{cases}
$$
If $k\ge0$ then $\lambda^{(k)}$ stands for the following partition of $n+k$ into $m+k$ parts:
$$
\lambda^{(k)}=(\ell_1,\ell_2,\dots,\ell_m,\underbrace{1,\dots,1}_{k\text{ times}})
$$
(in particular, $\lambda^{(0)}=\lambda$). If $\mu=(m_1,m_2,\dots,m_s)\in\Lambda_{r,s}$ then $W_{r,s,\mu}$ stands for the set of all words $\mathbf u$ such that:
\begin{itemize}
\item the length of $\mathbf u$ equals $r$;
\item $\mathbf u$ depends on the letters $x_1,x_2,\dots,x_s$;
\item for every $i=1,2,\dots,s$, the number of occurrences of $x_i$ in $\mathbf u$ equals $m_i$.
\end{itemize}
For a partition $\lambda=(\ell_1,\ell_2,\dots,\ell_m)\in\Lambda_{n,m}$, we put
$$
\mathbf S_\lambda=\var\bigl\{\mathbf{u\approx v}\mid\,\text{there is }i\in\bigl\{0,1,\dots,s(\lambda)\bigr\}\text{ such that }\mathbf u, \mathbf v\in W_{n+i,m+i,\lambda^{(i)}}\bigr\}.
$$
We call sets of the form $W_{n,m,\lambda}$ \emph{transversals}. We say that an overcommutative variety $\mathbf V$ \emph{reduces} [\emph{collapses}] a transversal $W_{n,m,\lambda}$ if $\mathbf V$ satisfies some non-trivial identity [all identities] of the form $\mathbf{u\approx v}$ with $\mathbf u,\mathbf v\in W_{n,m,\lambda}$. An overcommutative variety $\mathbf V$ is said to be \emph{greedy} if it collapses any transversal it reduces.

\begin{theorem}
\label{OC neutr,stand,costand,canc,distr,codistr}
For an overcommutative semigroup variety $\mathbf V$, the following are equivalent:
\begin{itemize}
\item[\textup{a)}]$\mathbf V$ is $\mathbb{OC}$-neutral;
\item[\textup{b)}]$\mathbf V$ is $\mathbb{OC}$-standard;
\item[\textup{c)}]$\mathbf V$ is $\mathbb{OC}$-costandard;
\item[\textup{d)}]$\mathbf V$ is $\mathbb{OC}$-cancellable;
\item[\textup{e)}]$\mathbf V$ is $\mathbb{OC}$-distributive;
\item[\textup{f)}]$\mathbf V$ is $\mathbb{OC}$-codistributive;
\item[\textup{g)}]$\mathbf V$ is greedy;
\item[\textup{h)}]either $\mathbf{V=SEM}$ or $\mathbf V=\mathop\bigwedge\limits_{i=1}^k\mathbf S_{\lambda_i}$ for some partitions $\lambda_1,\lambda_2,\dots,\lambda_k$.
\end{itemize}
\end{theorem}

The equivalence of the claims~a)--c) and~e)--g) of this theorem was proved in~\cite{Vernikov-01} (claim~g) was not mentioned in~\cite{Vernikov-01} explicitly but the fact that this claim is equivalent to each of the claims~a)--c),~e) and~f) readily follows from the proof of~\cite[Theorem~2]{Vernikov-01}). The equivalences~d)\,$\Longleftrightarrow$\,~g) and~g)\,$\Longleftrightarrow$\,~h) are verified in~\cite{Shaprynskii-Vernikov-21} and~\cite{Shaprynskii-Vernikov-11} respectively. The results of the paper~\cite{Volkov-94} and Propositions~\ref{Sub neutr,stand,costand,canc,distr,codistr} and~\ref{G-sets neutr,stand,costand,canc,distr,codistr} play the crucial role in the parts of the proof of Theorem~\ref{OC neutr,stand,costand,canc,distr,codistr} given in~\cite{Shaprynskii-Vernikov-21} and~\cite{Vernikov-01}. 

\subsection*{Acknowledgments.}\label{thanks} The author thanks Drs. Olga Sapir and Edmond W.\,H. Lee for a number of useful suggestions on improvement of the initial version of the survey.

\end{document}